\documentclass[aos,MSNbibl,seceqn,nameyear,dvips]{arximspdf}

\doi{10.1214/13-AOS1087} %kopijuoti is PTS
\volume{41}
\issue{1}
\pubyear{2013}
\firstpage{342}
\lastpage{369}

\makeatletter

\newproclaim{remark}{Remark}
\newproclaim{example}{Example}
\newcommand{\vu}{\mathbf{u}}
\newcommand{\vA}{\mathbf{A}}
\newcommand{\vB}{\mathbf{B}}
\newcommand{\vX}{\mathbf{X}}
\newcommand{\vY}{\mathbf{Y}}
\newcommand{\vD}{\mathbf{D}}

\newcommand{\vnull}{\mathbf{0}}

\newcommand{\vSigma}{\bolds{\Sigma}}
\newcommand{\vbeta}{\bolds{\beta}}
\newcommand{\vPi}{\bolds{\Pi}}

\newcommand{\E}{\mathrm{E}}

\newcommand{\ra}{\rightarrow}

\newproclaim{definition}{Definition}[section]

\newtheorem{lemma}[definition]{Lemma}
\newtheorem{theorem}[definition]{Theorem}

\newcommand{\lambdamax}{\lambda_{\mathrm{max}}}
\newcommand{\lambdamin}{\lambda_{\mathrm{min}}}
\newcommand{\defined}{\stackrel{\Delta}{=}}

\newcommand{\argmin}{\operatorname{arg\,min}}
\makeatother

\begin{document}
\begin{frontmatter}

\title{Quantile-adaptive model-free variable screening for
high-dimensional heterogeneous data}
\runtitle{Quantile-adaptive screening}

\begin{aug}
\author[A]{\fnms{Xuming} \snm{He}\thanksref{t2}\ead[label=e1]{xmhe@umich.edu}},
\author[B]{\fnms{Lan} \snm{Wang}\corref{}\thanksref{m1}\ead[label=e2]{wangx346@umn.edu}}
\and
\author[C]{\fnms{Hyokyoung Grace} \snm{Hong}\thanksref{t4}\ead[label=e3]{hhong@stt.msu.edu}}
\thankstext{t2}{Supported in part by NSF Grant DMS-10-07396, NIH Grant R01GM080503 and National Natural Science Foundation of China Grant 11129101.}
\thankstext{m1}{Supported in part by NSF Grant DMS-10-07603.}
\thankstext{t4}{Supported in part by a Baruch College Eugene M. Lang Fellowship.}

\runauthor{X. He, L. Wang and H. G. Hong}

\affiliation{University of Michigan, University of Minnesota
and Michigan~State~University}\address[A]{X. He\\
Department of Statistics \\
University of Michigan\\
Ann Arbor, Michigan 48109\\
USA\\
\printead{e1}}

\address[B]{L. Wang\\
School of Statistics\\
University of Minnesota\\
Minneapolis, Minnesota 55455\\
USA\\
\printead{e2}}

\address[C]{H. G. Hong\\
Department of Statistics and Probability \\
Michigan State University\\
East Lansing, Michigan 48824\\
USA\\
\printead{e3}}
\end{aug}

% HISTORY:
\received{\smonth{2} \syear{2012}}
\revised{\smonth{1} \syear{2013}}

% ABSTRACT
%
\begin{abstract}
We introduce a quantile-adaptive framework for nonlinear variable screening
with high-dimensional heterogeneous data. This framework has two
distinctive features:
(1) it allows the set of active variables
to vary across quantiles, thus making it more flexible to accommodate
heterogeneity; (2)
it is model-free
%in the sense of Zhu, Li, Li, Zhu (2011) as
and avoids the difficult
task of specifying the form of a statistical model in a high
dimensional space. Our nonlinear
independence screening procedure employs spline approximations to model
the marginal effects
at a quantile level of interest. Under appropriate conditions on the
quantile functions without
requiring the existence of any moments, the new procedure is shown to
enjoy the sure screening property in ultra-high dimensions.
Furthermore, the quantile-adaptive framework can naturally handle censored
data arising in survival analysis.
We prove that the sure screening property remains valid when the
response variable
is subject to random right censoring.
Numerical studies confirm the fine performance of the proposed method
for various semiparametric models and its effectiveness to extract
quantile-specific information
from heteroscedastic data.
\end{abstract}
%
% KEYWORDS

\begin{keyword}[class=AMS]
\kwd[Primary ]{68Q32}
\kwd{62G99}
\kwd[; secondary ]{62E99}
\kwd{62N99}
\end{keyword}

\begin{keyword}
\kwd{Feature screening}
\kwd{high dimension}
\kwd{polynomial splines}
\kwd{quantile regression}
\kwd{randomly censored data}
\kwd{sure independence screening}
\end{keyword}

\end{frontmatter}

%s1 #&#
\section{Introduction}\label{sec1}
We consider the problem of analyzing ultra-high dimensional data, where
the number of
candidate covariates (or features) may increase at an exponential rate.
Many efforts have been devoted to this challenging problem in recent
years, motivated
by modern applications in genomics, bioinformatics, chemometrics, among others.
A practically appealing approach is to first use a fast screening
procedure to reduce the dimensionality of
the feature space to a moderate scale; and then apply more
sophisticated variable selection techniques
in the second stage. In this paper, we propose a new quantile-adaptive,
model-free variable
screening procedure, which is particularly appealing for analyzing high
dimensional heterogeneous data and data with censored responses.

\citet{F2} proposed the sure independence screening (SIS)
methodology for linear regression
which screens variables by ranking their marginal correlations with the
response variable. They established the desirable sure screening
property, that is, some important features are retained with
probability approaching one, even if the dimensionality of the features
is allowed to grow exponentially fast
with the sample size.
\citet{F4} further extended the methodology to generalized
linear models, see also \citet{F3}. The problem of
nonlinear features screening
was addressed in \citet{Hall} using generalized correlation
ranking and more systematically
in \citet{F1} using nonparametric marginal ranking,
which extended the scope of applications of sure independence
screening. B\"{u}hlmann, Kalisch, and Maathuis (\citeyear{Buh-2010})
introduced the new concept of partial faithfulness and
proposed a computationally efficient
PC-algorithm for feature screening in linear models.

\citet{Zhu} proposed a novel feature screening procedure which
avoids the specification of a particular model structure. This
model-free screening framework is very appealing
because a misspecified model could easily corrupt the performance of a
variable selection method. Partly motivated by this interesting piece
of work, we propose a new framework
called quantile-adaptive model-free screening. We advocate a
quantile-adaptive approach which allows the
set of active variables to be different when modeling different
conditional quantiles.
This new framework provides a more complete picture of the conditional
distribution of the response
given all candidate covariates and is more natural and effective for
analyzing high-dimensional
data that are characterized by heteroscedasticity.

In the quantile-adaptive model-free screening framework, we estimate
marginal quantile regression
nonparametrically using
$B$-spline approximation.
In this aspect, our technique shares some similarity with
that in \citet{F1}. The main technical challenge
is to deal with the nonsmooth loss function, because the nonparametric
marginal utility
we consider does not have a closed form expression as in \citet{F1}. We derive useful exponential
bounds using the empirical process theory to establish the sure
screening property.
When working with marginal quantile regression,
the usual sub-Gaussian tail type condition in high-dimensional analysis
can be relaxed
and replaced by the assumption that
the conditional density function of the random error
has a positive lower bound around the quantile of interest.
Empirically, we also demonstrate that
the proposed procedure works well with heavy-tailed error distributions.

Sure independence screening remains challenging and little explored
when the response variable
is subject to random censoring, a common problem in survival analysis.
\citet{F12} extended
the methodology of sure independence screening using the marginal
Cox proportional hazards model and studied its performance empirically.
In this paper,
we demonstrate that in the quantile-adaptive model-free screening
framework, randomly censored
responses can be naturally accommodated by ranking a marginal weighted
quantile regression
utility. We establish the sure screening property for the censored case under
some general conditions.

The rest of the paper is organized as follows. In Section~\ref{sec2}, we
introduce the quantile-adaptive model-free feature screening procedure.
In Section~\ref{sec3}, we investigate its theoretical properties. Section~\ref{sec4}
discusses the extension
to survival analysis. In Section~\ref{sec5}, we carry out simulation studies to access
the performance of the proposed method. The numerical results
demonstrate the favorable performance of the proposed method,
especially when the errors are heavy-tailed or heteroscedastic. In
Section~\ref{sec6},
we demonstrate the application on a real data example. Section~\ref{sec7}
contains further discussions.
The technical details are given in Section~\ref{sec8}.

%s2 #&#
\section{Quantile-adaptive model-free feature screening}\label{sec2}

%s2.1 #&#
\subsection{A general framework}\label{sec2.1}

We consider the problem of nonlinear variable screening in
high-dimensional feature space, where we observe
a response variable $Y$ and associated covariates $X_1,\ldots, X_p$.
The goal is to rapidly
reduce the dimension of the covariate space $p$ to a moderate scale
via a computationally convenient procedure.
Since ultra-high dimensional data often display heterogeneity,
we advocate a quantile-adaptive feature screening framework. More
specifically, we assume
that at each quantile level a sparse set of covariates are relevant for
modeling $Y$, but allow this set to be different at different
quantiles, see, for instance, Examples~\ref{ex2} and~\ref{ex3} in Section~\ref{sec5}.
At a given quantile level $\alpha$ ($0<\alpha<1$), we define the set
of active variables
\[
M_{\alpha}=\bigl\{j\dvtx Q_{\alpha}(Y|{\vX}) \mbox{ functionally depends on } X_j \bigr\},
\]
where $Q_{\alpha}(Y|\vX)$ is the $\tau$th conditional quantile of
$Y$ given $\vX=(X_1,\ldots,X_p)^T$, that is,
$Q_{\alpha}(Y|\vX)=\inf\{y\dvtx P(Y\leq y|\vX)\geq\alpha\}$.
Let $S_{\alpha}=|M_{\alpha}|$ be the cardinality of~$M_{\alpha}$.
Throughout this paper, we assume $S_{\alpha}$,
$0<\alpha<1$, is smaller than the sample size $n$.

%To illustrate the definition of active variables, we consider
%the following semiparametric heteroscedastic regression model
%Y=\valpha_1^T\vX_{A_1}+f_1(\valpha_2^T\vX_{A_2})+f_2(\valpha_3^T
%where $\valpha_i$, $i=1,2,3$, are unknown vectors of coefficients;
%$f_i$, $i=1,2$, are unknown smooth functions;
%$A_i$ are subsets of $\{1,2,\ldots,p\}$; $\vX_{A_i}$
%denotes the subvector of $\vX$ consisting of those entries with
%indices in $A_i$;
%and $\varepsilon$ has a distribution symmetric about 0.
%Then at
%the median ($\tau=0.5$), the set of active variables comprises the
%covariates with indices in $A_1\cup A_2$;
%while at the first quartile ($\tau=0.25$),
%it comprises the covariates with indices in $A_1\cup A_2 \cup A_3$. }

In practice, we may consider several quantiles to explore the sparsity
pattern and the effects of the covariates
at different parts of the conditional distribution.
%In contrast, existing screening procedures only focus on
%the center of the conditional distribution and may ignore
%heterogeneity which is often of direct importance.
We refer to \citet{K2} for a comprehensive introduction to quantile
regression.
%and to He (2009) for a recent general overview.

In ultrahigh dimensional data analysis,
there generally exists little prior information for specifying
a\vadjust{\goodbreak}
statistical model.
Given a large number of covariates, it is hard to determine which
covariates have linear effects
and which have nonlinear effects. In our framework besides the sparsity
assumption, we do not impose a specific model structure but allow the
covariate effects to be nonlinear.

%Since
%independent}\},

%s2.2 #&#
\subsection{Ranking by marginal quantile utility}\label{sec2.2}
Let $\{(\vX_i, Y_i), i=1,\ldots, n\}$ be i.i.d.
copies of $(\vX, Y)$, where $\vX_i=(X_{i1},\ldots,X_{ip})^T$. Note that
\begin{eqnarray*}
Y \mbox{ and } X_j \mbox{ are independent}\quad\Leftrightarrow\quad
Q_{\alpha}(Y|X_j)-Q_{\alpha}(Y)=0\qquad \forall\alpha
\in(0,1),
\end{eqnarray*}
where $Q_{\alpha}(Y|X_j)$ is the $\alpha$th conditional quantile of
$Y$ given $X_j$ and
$Q_{\alpha}(Y)$ is the $\alpha$th unconditional quantile of $Y$.
To estimate the effect of $X_j$ on~$Y$, we consider the marginal quantile
regression of $Y$ on $X_j$.
Let $f_j(X_j)=\argmin_{f}\E[\rho_{\alpha}(Y-f(X_j))-\rho_\alpha(Y )]$, where the inclusion of $\rho_\alpha(Y )$ makes the expectation
well defined even when $Y$ has no finite moment, where
$\rho_{\alpha}(u) = u \{\alpha- I(u<0) \}$ is the quantile
loss function (or check function).
It is known that $f_j(X_j)=Q_{\alpha}(Y|X_j)$, the $\alpha$th
conditional quantile of $Y$ given $X_j$.

Without loss of generality, we assume that each $X_{j}$ takes values on
the interval $[0,1]$.
Let $\mathbb{F}$ be the class of functions defined in condition (C1)
in Section~\ref{sec3.1}.
Let
$0=s_0<s_1<\cdots<s_k=1$ be a partition of the interval. Using the $s_i$
as knots, we construct $N=k+l$ normalized $B$-spline basis functions of
order $l+1$
which form a basis for $\mathbb{F}$. We write these basis
functions as a vector $\bolds{\pi}(t)=(B_1(t),\ldots
,B_N(t))^T$, where $\Vert B_k(\cdot)\Vert_{\infty}\leq1$ and
$\Vert \cdot\Vert_{\infty}$ denotes the sup norm.
Assume that $f_j(t)\in\mathbb{F}$. Then $f_j(t)$ can be well approximated
by a linear combination of the basis functions $\bolds{\pi
}(t)^T\vbeta$, for some $\vbeta\in\mathbb{R}^N$.

Let $
\widehat{\vbeta}_j=\argmin_{\vbeta\in\mathbb{R}^N}\sum_{i=1}^n\rho_{\alpha}(Y_i-\bolds{\pi}(X_{ij})^T\vbeta),
$ and define
\[
\widehat{f}_{nj}(t)=\bolds{\pi}(t)^T\widehat{\vbeta
}_{j}-F_{Y,n}^{-1}(\alpha),
\]
where $F_{Y,n}^{-1}(\alpha)$ is the $\alpha$th sample quantile
function based on $Y_1,\ldots,Y_n$.
Thus $\widehat{f}_{nj}(t)$
is a nonparametric estimator of $Q_{\alpha}(Y|X_j)-Q_{\alpha}(Y)$.
We expect $\widehat{f}_{nj}$ to be close to zero if $X_j$ is
independent of $Y$.
%Note that $\widehat{f}_{nj}(t)$ is the sample version of
%$f_{nj}(t)$.

The independence screening is based on the magnitude of the estimated
marginal components
$\Vert \widehat{f}_{nj}\Vert_n^2=n^{-1}\sum_{i=1}^n\widehat{f}_{nj}(X_{ij})^2$.
More specifically, we will select the subset of variables
\[
\widehat{M}_{\alpha}=\bigl\{1\leq j\leq p\dvtx \Vert \widehat{f}_{nj}
\Vert_n^2\geq \nu_n\bigr\},
\]
where $\nu_n$ is a predefined threshold value. In practice, we often
rank the features
by $\Vert \widehat{f}_{nj}\Vert_n^2$ and keep the top $[n/\log(n)]$
features, where
$[a]$ denotes the integer part of~$a$.

%s3 #&#
\section{Theoretical properties}\label{sec3}

%s3.1 #&#
\subsection{Preliminaries}\label{sec3.1}
We impose the following regularity conditions
to facilitate our technical derivations.\vadjust{\goodbreak}
\begin{longlist}[(C1)]
\item[(C1)] The conditional quantile function $Q_{\alpha}(Y|X_j)$ belongs to
$\mathbb{F}$, where
$\mathbb{F}$ is the class of functions defined on $[0,1]$ whose $l$th
derivative satisfies a Lipschitz condition of order $c$:
$|f^{(l)}(s)-f^{(l)}(t)|\leq c_0|s - t|^c$, for some positive constant $c_0$, $s,t\in[0,1]$, where $l$ is a nonnegative integer and $c\in(0,1]$
satisfies $d=l+c>0.5$.

\item[(C2)] $\min_{j\in M_{\alpha}}E{(Q_{\alpha}(Y|X_j)-Q_{\alpha
}(Y))^2}\geq c_1n^{-\tau}$ for some $0\leq\tau<\frac{2d}{2d+1}$ and
some positive constant $c_1$.

%We assume that $\min_{j\in M_{*}}E{f_j(X_j)^2}\geq c_1Nn^{-2\kappa}$
%for some $0<\kappa<\frac{d}{2d+1}$
%and some positive constant $c_1$.\\
%(3) Assume that $\varepsilon_i$ satisfies condition $E(\exp(M_1|
%for some positive constants $M_1$ and $M_2$.\\
\item[(C3)] The conditional density $f_{Y|X_j}(t)$ is bounded away from 0 and
$\infty$ on $[Q_{\alpha}(Y|X_j)-\xi, Q_{\alpha}(Y|X_j)+\xi]$, for
some $\xi>0$, uniformly in $X_j$.

\item[(C4)] The marginal density function $g_j$ of $X_j$, $1\leq j \leq p$,
are uniformly bounded away from 0 and $\infty$.

\item[(C5)] The number of basis functions $N$ satisfies $N^{-d}n^{\tau}=o(1)$ and
$Nn^{2\tau-1}=o(1)$ as $n\ra\infty$.
%$0<K_1\leq g_j(X_j)\leq K_2<\infty$ on $[a,b]$ for $1\leq j \leq p$
%for some constants $K_1$ and $K_2$.}
\end{longlist}

Condition (C1) assumes that the conditional quantile function
$Q_{\alpha}(Y|X_j)$ belongs to a class of smooth functions. This
condition is standard for nonparametric spline approximation. Condition
(C2) assumes that the features in the active set at quantile level
$\alpha$ have strong enough marginal signals;
a smaller $\tau$ corresponds to a stronger marginal signal.
This condition is important as it guarantees that marginal utilities
carries information about the features in
the active set.
Condition (C3) is a standard condition on random errors in the theory
for quantile regression.
It relaxes the usual sub-Gaussian assumptions that are needed in the
literature on high dimensional inference.
%Note that in high-dimensional inference it usually needs to assume
%that the random errors follow the sub-Gaussian
%distribution which is essential for deriving the exponential
%probability bounds. \bfblue{Condition} (C3)
%relaxes this assumption and allows heavy-tailed error distributions.
Condition (C4) is similar as condition (B) of \citet{F1}.
Note that (C4) is not restrictive when $X_j$ is supported on a bounded
interval, say $[0,1]$. When the distribution of $X_j$ has an unbounded
support (e.g., normal), we
can view $X_j$ as coming from a truncated distribution. In fact, if
there is an outlier in $X_j$ as in the case of a heavy-tailed
distribution, we do better by dropping the outlier or transforming
$X_j$ to be uniformly distributed on $[0,1]$.
In our numerical simulations of Section~\ref{sec5}, the normally distributed
covariates are scaled to the interval $[0,1]$ and the results would
change little if any reasonable truncation is used instead.
Condition (C5) describes how fast the number of basis functions is allowed
to grow with the sample size.

Given $(Y,\vX)$, where $\vX=(X_1,\ldots,X_p)^T$,
we define
%
%e3.1 #&#
\begin{equation}
\label{beta0} \vbeta_{0j}=\mathop{\argmin}_{\vbeta\in\mathbb{R}^N}\E\bigl[
\rho_{\alpha
}\bigl(Y-\bolds{\pi}(X_{j})^T\vbeta\bigr)
-\rho_\alpha(Y )\bigr].
\end{equation}
Let $f_{nj}(t)=\bolds{\pi}(t)^T\vbeta_{0j}-Q_{\alpha
}(Y)$, whose sample version was defined in Section~\ref{sec2}.
Furthermore, we let $\Vert f_{nj}\Vert^2=E[f_{nj}(X_j)^2]$. The following
lemma shows that the spline
approximation error is negligible in the sense that the spline
approximation carries the
same level of information about the marginal signal.

%We assume that $\min_{j\in M_{\alpha}}E{(Q_{\alpha}(Y|X_j)-Q_{
%and some positive constant $c_1$.

%le3.1 #&#
\begin{lemma}\label{lemma1}
%Assume that $N^{-2d}n^{\tau}=o(1)$. Then
Under condition \textup{(C5)},
$\min_{j\in M_{\alpha}}\Vert f_{nj}\Vert^2\geq c_1n^{-\tau}/8$,
for all $n$ sufficiently large.
\end{lemma}

%s3.2 #&#
\subsection{Sure screening property}\label{sec3.2}
As covariate screening often serves only as the first step for high
dimensional data analysis,
the most important property as far as practical application is
concerned is the sure screening property.
In the quantile-adaptive framework, we require that the sure screening
property holds at each quantile
level $\alpha$, that is, the set of selected covariates at quantile
level $\alpha$ includes $M_{\alpha}$ with probability tending to one.

The key step of deriving the sure screening\vspace*{1.5pt} property is to establish
exponential probability bounds
for $\Vert \widehat{\vbeta}_j-\vbeta_{0j}\Vert $ and $\Vert \widehat
{f}_{nj}\Vert_n^2-\Vert f_{nj}\Vert^2$. The main technical challenge is that
$\widehat{\vbeta}_j$
is defined by minimizing a nonsmooth objective function, thus does not
have a closed-form
expression. The exponential bounds are summarized in the following lemma.

%le3.2 #&#
\begin{lemma}\label{lemma0}
Assume conditions \textup{(C1)--(C5)} are satisfied.
\begin{longlist}[(1)]
%$n^{2\tau-1}\log n=o(1)$
%in the original statement can be dropped.)}\\
\item[(1)]
For any $C>0$, there exist positive constants $c_2$ and $c_3$ such that
\begin{eqnarray*}
&&P\Bigl(\max_{1\leq j\leq p}\bigl\Vert \widehat{\vbeta}_j-
\vbeta_{0j}\bigr\Vert \geq CN^{1/2}n^{-\tau}\Bigr)
\\
&&\qquad \leq 2p\exp \bigl(-c_2n^{1-4\tau} \bigr) +p\exp
\bigl(-c_3N^{-2}n^{1-2\tau} \bigr)
\end{eqnarray*}
for all $n$ sufficiently large.

\item[(2)] For any $C>0$, there exist some positive constants $\delta_1$ and
$\delta_2$
such that
\begin{eqnarray*}
&&P \Bigl(\max_{1\leq j\leq p} \bigl| \bigl\Vert \widehat {f}_{nj}
\bigr\Vert_n^2-\Vert f_{nj}\Vert^2 \bigr| \geq
Cn^{-\tau} \Bigr)
\\
&&\qquad \leq p \bigl\{11\exp \bigl(-\delta_1n^{1-4\tau} \bigr)
+12N^2\exp \bigl(-\delta_2N^{-3}n^{1-2\tau}
\bigr) \bigr\},
\end{eqnarray*}
for all $n$ sufficiently large.
\end{longlist}
\end{lemma}

%imposing extra condition on $\tau$. For example, if as in the second
%part of
%Theorem 3.3, we further assume $\tau<1/4$, then
%the constant 11 can be absorbed into $\delta_1$. In the paper of Fan
%et al, similar constants also appear
%in their bounds.)}

\begin{remark*} The results suggest that we can handle the
dimensionality $\log p=o(n^{1-4\tau}+N^{-3}n^{1-2\tau})$. This dimensionality
depends on the number of basis functions $N$ and the strength of the
marginal signals.
If we take $N=n^{1/(2d+1)}$ (the optimal rate
for spline approximation), then for $\tau<\min(1/4,(d-1)/(2d+1))$, we
can handle
ultra-high dimensionality, that is, $p$ can grow at the exponential rate.
\end{remark*}

The following theorem establishes the sure screening property.

%th3.3 #&#
\begin{theorem}[(Sure screening property)]\label{main1}
%Assume conditions (C1)-(C4), $N^{-d}n^{\tau}=o(1)$,
%$Nn^{2\tau-1}=o(1)$ and $n^{2\tau-1}\log n=o(1)$.
Under the conditions of Lem\-ma~\ref{lemma0},
if $\tau<1/4$, $N^{3}n^{2\tau-1}=o(1)$, and
we take the threshold value $\nu_n=\delta^* n^{-\tau}$ with $\delta^*\leq c_1/16$
for the constant $c_1$ specified in Lemma~\ref{lemma1}, then
\begin{eqnarray*}
P (M_{\alpha}\subset\widehat{M}_{\alpha} ) \geq1-S_{\alpha}
\bigl\{11\exp \bigl(-\delta_1n^{1-4\tau} \bigr)+12N^2
\exp \bigl(-\delta_2N^{-3}n^{1-2\tau} \bigr) \bigr\},
\end{eqnarray*}
for all $n$ sufficiently large. Especially, $P (M_{\alpha}\subset
\widehat{M}_{\alpha} ) \ra1$
as $n\ra\infty$.
\end{theorem}

%s3.3 #&#
\subsection{Controlling false discovery}\label{sec3.3}

An interesting question is how many variables are retained after the
screening. A simple bound is provided below, which extends the results
in \citet{F1}.

Let $\vPi= (\bolds{\pi}(X_1), \ldots,\bolds{\pi}(X_p)
)^T$ and $\vSigma=\E (\vPi\vPi^T )$.
Let $\Vert \cdot\Vert_{F}$ denote the Frobenius matrix norm
and let $\Vert \cdot\Vert_2$ denote the spectral matrix norm.
Note that
\begin{eqnarray*}
\sum_{j=1}^p\Vert f_{nj}
\Vert^2&=& \sum_{j=1}^p\E
\bigl(\bolds{\pi }(X_j)\vbeta_{0j} \bigr)^2
\leq\sum_{j=1}^p \lambdamax \bigl(\E\bolds{
\pi}(X_j)\bolds{\pi }(X_j)^T \bigr)\Vert
\vbeta_{0j}\Vert^2
\\
&\leq& N\sum_{j=1}^p \mbox{trace} \bigl(
\E\bolds{\pi}(X_j)\bolds {\pi}(X_j)^T \bigr)
\leq N\E \Biggl[\sum_{j=1}^p\sum
_{k=1}^NB_{k}^2(X_j)
\Biggr]
\\
&\leq& N\E \bigl(\Vert \mathbf{\vPi}\Vert_{F}^2 \bigr)
\leq N^2 \E \bigl(\Vert \mathbf{\vPi}\Vert_{2}^2
\bigr)=N^2\lambdamax (\vSigma ),
\end{eqnarray*}
where the second inequality uses the result
$\Vert \vbeta_{0j}\Vert^2\leq cN$ for some positive constant $c$ (proved in
the supplemental material [\citet{HWH}]).

For any $\varepsilon>0$, we define the set
\[
D_n= \Bigl\{\max_{1\leq j\leq p} \bigl| \Vert \widehat {f}_{nj}
\Vert_n^2-\Vert f_{nj}\Vert^2 \bigr|\leq
\varepsilon n^{-\tau} \Bigr\}.
\]
Then on $D_n$,
the cardinality of $\{j\dvtx \Vert \widehat{f}_{nj}\Vert^2_n>2\varepsilon n^{-\tau
}\}$
cannot exceed the cardinality of $\{j\dvtx \Vert f_{nj}\Vert^2>\varepsilon n^{-\tau
}\}$, which is bounded by
$\varepsilon^{-1}N^2n^{\tau}\lambdamax (\vSigma )$.

For
$\nu_n=\delta^* n^{-\tau}$, we take
$\varepsilon=\delta^*/2$, then
$P (|\widehat{M}_{\alpha}|\leq\varepsilon^{-1} N^2n^{\tau
}\lambdamax (\vSigma ) )
\geq P(D_n).
$ Applying Lemma~\ref{lemma0}(2), we obtain the following theorem which provides
a bound on the size of selected variables.

%th3.4 #&#
\begin{theorem}\label{bound1}
%Assume conditions (C1)-(C4), $N^{-d}n^{\tau}=o(1)$,
%$Nn^{2\tau-1}=o(1)$ and $n^{2\tau-1}\log n=o(1)$. Then for
%$\nu_n=\delta^* n^{-\tau}$ with $\delta^*\leq c_1/16$,
Under the conditions of Theorem~\ref{main1}, there
exist some positive constants $\delta_1$ and $\delta_2$ such that
for all $n$ sufficiently large,
\begin{eqnarray*}
&& P \bigl(|\widehat{M}_{\alpha}|\leq2N^2n^{\tau}
\lambdamax (\vSigma )/\delta^* \bigr)
\\
&&\qquad \geq 1- p \bigl\{11\exp \bigl(-\delta_1n^{1-4\tau} \bigr)
+12N^2\exp \bigl(-\delta_2N^{-3}n^{1-2\tau}
\bigr) \bigr\}.
\end{eqnarray*}
Especially, $P (|\widehat{M}_{\alpha}|\leq2N^2n^{\tau
}\lambdamax (\vSigma )/\delta^* )\ra1$
as $n\ra\infty$.
\end{theorem}

The above theorem suggests that if $\lambdamax (\vSigma
)=O(n^{\gamma})$ for some $\gamma>0$, then the model obtained after
screening is of polynomial size with high probability.
%It thus effectively reduces the dimensionality from the exponential
%size.
Similar observation has been reported for the $L_2$-based screening
procedure of Fan, Feng and Song (\citeyear{F1}).

%s4 #&#
\section{Quantile-adaptive screening in survival analysis}\label{sec4}
There exists very limited amount of work on feature screening with
censored responses.
\citet{F12} and Zhao and Li (\citeyear{Zhao}) investigated marginal
screening based on
the Cox proportional hazards model. As a powerful alternative to the\vadjust{\goodbreak}
classical Cox model,
quantile regression has recently emerged as a useful tool for analyzing
censored data, see
\citet{Ying-1995}, McKeague, Subramanian and Sun (\citeyear{McKeague-2001}),
\citet{Portnoy-2003}, \citet{peng}, \citet{wang}, among others.
The quantile regression approach directly models the survival time and
is easy to interpret.
Furthermore, it relaxes the
proportional hazards assumption of the Cox model and can naturally
accommodate heterogeneity in the data.
The quantile regression based screening procedure can be naturally
extended to survival analysis.

%Fan, Feng and Wu \bfblue{(2010)} proposed marginal Cox model based
%screening procedure but the theory is
%still lacking.
Assume that $Y_i$ is subject to random right censoring. Instead
of $\{(\vX_i, Y_i),\break i=1,\ldots, n\}$, we observe $\{(\vX_i, Y_i^*,
\delta_i), i=1,\ldots, n\}$
where
%
%e4.1 #&#
\begin{equation}
Y_i^*=\min(Y_i,C_i),\qquad
\delta_i=I(Y_i\leq C_i). \label{surv}
\end{equation}
The random variable $C_i$, called the censoring variable, is assumed to
be conditionally
independent of $Y_i$ given $\vX_i$. In this section, we assume that
the censoring distribution is the same for all covariates, but this
assumption will be relaxed in Section~\ref{sec7.1}. Let $G(t)=P(C_i>t)$ be
the survival function of $C_i$. Let $\widehat{G}(t)$ be the
Kaplan--Meier estimator of $G(t)$,
based on $\{Y_i^*, \delta_i\}$, $i=1,\ldots,n$.

%(The above two conditions imply that $P(C\leq c, Y\leq y|\vX)=P(C\leq
%c)P(Y\leq y|\vX)$. It seems to be weaker
%than the condition $C$ being independent of $(\vX, Y)$, as the latter
%implies the former but the former does not
%necessarily imply the latter. Or is it not true?)
%}

Similarly as in the case of complete data, we
consider independence screening based on nonparametric marginal
regression given $X_j$.
More specifically, we consider inverse probability weighted marginal
quantile regression estimator
\[
\widehat{\vbeta}_j^{c}=\mathop{\argmin}_{\vbeta\in\mathbb{R}^N}\sum
_{i=1}^n\frac{\delta_i}{\widehat{G}(Y_i^*)}
\rho_{\alpha
}\bigl(Y_i^{*}-\bolds{
\pi}(X_{ij})^T\vbeta\bigr).
\]
Let
$
\widehat{f}_{nj}^{c}(t)=\bolds{\pi}(t)^T\widehat{\vbeta
}_{j}^{c}-F_{\mathrm{KM},n}^{-1}(\alpha)
$
where $F_{\mathrm{KM},n}^{-1}(\alpha)$ is the nonparametric estimator
of the $\alpha$th
conditional quantile of $Y$ based on $(Y_i^*, C_i, \delta_i)$,
$i=1,\ldots,n$. The estimator
we use here is the inverse function (the left-continuous version) of
the Kaplan--Meier estimator of the distribution function of
$Y$, whose properties have been studied in \citet{Lo}. We will
select the subset of variables
$
\widehat{M}_{\alpha}^{c}=\{1\leq j\leq p\dvtx \Vert \widehat
{f}_{nj}^{c}\Vert_n^2\geq\nu_n^{c}\},
$
where $\nu_n^{c}$ is a predefined threshold value. As for the complete
data case,
in practice we often rank the features
by $\Vert \widehat{f}^{c}_{nj}\Vert_n^2$ and keep the top $[n/\log(n)]$ features.

For the random censoring case, in addition to conditions (C1)--(C5), we
assume that:
\begin{longlist}[(C5)]
\item[(C6)] $P(t\leq Y_i\leq C_i)\geq\tau_0>0$ for some positive constant
$\tau_0$
and any $t\in[0,T]$, where $T$ denotes the maximum follow-up time.
Furthermore, $\sup\{t\dvtx\break P(Y>t)>0\}\geq\sup\{t\dvtx P(C>t)>0\}$.
The survival function of the censoring variable $G(t)$ has uniformly
bounded first derivative.

\item[(C7)] There exist $0< \beta_1< \beta_2 <1$ such that $\alpha\in
[\beta_1 , \beta_2]$ and that the distribution function of $Y_i$ is
twice differentiable in $[Q_{\beta_1}(Y_i)-\varepsilon, Q_{\beta
_2}(Y_i)+\varepsilon]$
for some $0<\varepsilon<1$, with the first derivative bounded away from zero
and the second derivative bounded in absolute value.
\end{longlist}

Condition (C6) and (C7) are common in the survival analysis literature
to ensure that
the Kaplan--Meier estimator and its inverse function are well behaved.
The following theorem states that
the sure screening property holds for the random censoring case and an upper
bound that controls the size of selected variables can be obtained.

%th4.1 #&#
\begin{theorem}\label{main2}
Assume conditions \textup{(C1)--(C7)} are satisfied, $\tau<1/4$ and
$N^{3}n^{2\tau-1}=o(1)$,
%$N^2n^{2\tau-1}\log n=o(1)$ and $n^{2\tau-3/4}(\log n)^{3/4}=o(1)$.
if we take the threshold value $\nu_n^c=\delta^* n^{-\tau}$ with
$\delta^*\leq c_1/16$, then:

\begin{longlist}[(1)]
\item[(1)] there exist positive constants $\delta_3$ and $\delta_4$ such that
\begin{eqnarray*}
P \bigl(M_{\alpha}\subset\widehat{M}_{\alpha}^{c} \bigr)
\geq1-S_{\alpha} \bigl\{17\exp \bigl(-\delta_3n^{1-4\tau}
\bigr) +12N^2\exp \bigl(-\delta_{4}N^{-3}n^{1-2\tau}
\bigr) \bigr\},
\end{eqnarray*}
for all $n$ sufficiently large. Especially,
$
P (M_{\alpha}\subset\widehat{M}_{\alpha}^{c} )\ra1
$ as $n\ra\infty$.

\item[(2)] for all $n$ sufficiently large,
\begin{eqnarray*}
&& P \bigl(\bigl|\widehat{M}_{\alpha}^{c}\bigr|\leq2N^2n^{\tau}
\lambdamax (\vSigma )/\delta^*  \bigr)
\\
&&\qquad \geq 1-p \bigl\{17\exp \bigl(-b_7n^{1-4\tau} \bigr)
+12N^2\exp \bigl(-b_{8}N^{-3}n^{1-2\tau}
\bigr) \bigr\}.
\end{eqnarray*}
Especially,
$P (|\widehat{M}_{\alpha}^{c}|\leq2N^2n^{\tau}\lambdamax(\vSigma )
/\delta^*   )\ra1$
as $n\ra\infty$.
\end{longlist}
\end{theorem}

%%%%%%%%%%%%%%%%%%%%%%%%%%
%s5 #&#
\section{Monte Carlo studies}\label{sec5}
We carry out simulation studies to investigate the performance
of the proposed quantile adaptive sure independence screening
procedure (to be denoted by QaSIS).
We consider two criteria for evaluating the performance as in \citet{Zhu}.
The first criterion is the minimum model size (denoted by
$\mathcal{R}$), that is, the smallest number of covariates that we
need to include to ensure that all the active variables are selected.
The second criterion is the proportion of active variables (denoted by
$\mathcal{S}$)
selected by the screening procedure when the threshold $\nu_n=[n/\log
(n)]$ is adopted. Note that
the first criterion does not need to specify a threshold. An effective variable
screening procedure is expected to have the value of $\mathcal{R}$
reasonably small comparing
to the number of active variables
and the value of $\mathcal{S}$ close to one.

We first consider the complete data case and compare the performance of QaSIS
with the nonparametric independence screening (NIS) procedure
of Fan, Feng and Song (\citeyear{F1}) and the sure independent ranking and
screening (SIRS) procedure
of \citet{Zhu}. In computing QaSIS and NIS, the number of basis
($d_n$) is set to be $[n^{1/5}] =3$. For each example, we report
the results based on 500 simulation runs.

\begin{example}[(Additive model, $n=400$, $p=1000$)]\label{ex1}
This example is
adapted from \citet{F1}. Let $g_1(x)=x$,
$g_2(x)=(2x-1)^2$, $g_3(x)={\sin(2 \pi x)}/{(2-\sin(2 \pi x))}$, and
$g_4(x)=0.1\sin(2 \pi x)+0.2 \cos(2 \pi x)+\break 0.3 \sin(2 \pi
x)^2+0.4 \cos(2 \pi x)^3+ 0.5 \sin(2 \pi x)^3$.
The following cases are studied:
\begin{itemize}
\item \textit{Case} (1a): $Y=5 g_1(X_1)+3g_2(X_2)+4g_3(X_3)+6g_4(X_4)+\sqrt{1.74}
\varepsilon$,
where the vector of covariates
$\mathbf{X}=(X_1,\ldots,X_{1000})^T$ is generated from the
multivariate normal distribution
with mean $\vnull$ and the covariance matrix
$\bolds{\Sigma}=(\sigma_{ij})_{1000 \times1000}$ with $\sigma_{ii}=1$ and
$\sigma_{ij}=\rho^{|i-j|}$ for $i\neq j$,
$\varepsilon\sim N(0,1)$ is independent of $\mathbf{X}$. In case (1a), we consider
$\rho=0$.
\item \textit{Case} (1b): same as case (1a) except that $\rho=0.8$.
\item \textit{Case} (1c): Same as case (1b) except that $\varepsilon$ has the Cauchy
distribution.
\end{itemize}
Note that the models are homoscedastic in Example~\ref{ex1}, thus the number of
active variables are the same across different quantiles.
\end{example}

\begin{example}[(Index model, $n=200$, $p=2000$)]\label{ex2}
This example is
adapted from
\citet{Zhu}. The random data are generated from
$Y=2(X_1+0.8X_2+0.6X_3+0.4X_4+0.2X_5)+\exp(X_{20}+X_{21}+X_{22}) \cdot
\varepsilon,
$
where $\varepsilon\sim N(0,1)$, $\vX=(X_1,X_2,\ldots,X_{2000})^T$
follows the
multivariate normal distribution with the correlation structure described
in case (1b).
%Cauchy distribution.
%The sample size is $n$=200 and number of candidate covariates is
%$p$=2000.
Different from the regression models in Example~\ref{ex1}, this model is
heteroscedastic:
the number of active variables is 5 at the median but 8 elsewhere.
\end{example}

\begin{example}[(A more complex structure, $n=400$, $p=5000$)]\label{ex3}
We
consider a more complex
heteroscedastic model for which the conditional distribution of $Y$
does not have a simple
additive or index structure.

\begin{itemize}
\item \textit{Case} (3a):
$Y=2(X_1^2+X_2^2)+
\{10^{-1}\exp(X_1+X_2+X_{18}+X_{19}+
\cdots+X_{30})\} \cdot\varepsilon$, where $\varepsilon\sim N(0,1)$,
and $\vX=(X_1,X_2,\ldots,X_{5000})^T$ follows the
multivariate normal distribution with the correlation structure described
in case (1b). In
this case, the number of active variables is 2 at the median but is 15
elsewhere.
\item \textit{Case} (3b): same as case (3a), but with
$2(X_1^2 + X_2^2)$ replaced by $2((X_1 + 1)^2 + (X_2 + 2)^2)$.
%we generate random data from $Y=2(X_1^2+X_2^2)+
%$\varepsilon(\alpha)=\varepsilon-F_{\chi^2(1)} (\alpha)$ with $\varepsilon\sim
%$\vX=(X_1,X_2,\ldots,X_{2000})^T$ follows the same distribution as in
%Case 3a. In
%Case 3b, the number of active variables is 2 at the $\alpha$th
%quantile but 15 elsewhere.
\end{itemize}
\end{example}

\begin{table}
\caption{Results for Examples \protect\ref{ex1}--\protect\ref{ex3}.
The numbers reported are the median of $\mathcal{R}$ [with
interquartile range (IQR) given in parentheses] and $\mathcal{S}$}\label{tab1}
\begin{tabular*}{\textwidth}{@{\extracolsep{\fill}}lccccc@{}}
%},
%label=t1, pos=h!, ]{lclclc}
%{
%{$p^{*}$: the number of truly active variables.}
%%{In Example 2, NIS(L)}\\ {is for the location part only and NIS (LS)
%%is for the location-scale part.}
%%\tnote[2]{In Example 2, NIS(L)}\\ {is for the location part only and
%%NIS (LS) is for the location-scale part.}
%%\tnote[3]{\red{Example 1--3 use the different $n$ and $p$: Example 1
%%($n=400, p=1000$), Example 2 ($n=200$, %$p=2000$), Example 3 ($n=400$,
%%$p=2000$).}}
%}
%{
\hline
\textbf{Example}&
\textbf{Case}&
\textbf{Method}&
\multicolumn{1}{c}{$\bolds{p^*}$}&
\multicolumn{1}{c}{$\bolds{\mathcal{R}}$\textbf{(IQR)}}&
\multicolumn{1}{c@{}}{$\bolds{\mathcal{S}}$}\\
\hline
Example~\ref{ex1} &(1a) &QaSIS $({\alpha=0.50})$ & \phantom{0}4&655 (434) & 0.56\\
&&QaSIS $({\alpha=0.75})$ &\phantom{0}4& 652 (398) & 0.56\\
&&NIS & \phantom{0}4&660 (415) & 0.55\\
&&SIRS &\phantom{0}4& 689 (365) & 0.53\\[3pt]
%%%%%%%%%%%%%%%%%%%%%%%%%%%
%1b &QaSIS$({\alpha=.50})$&4&$ 37.0~ (119.3)$ & 0.91 \NN
% &QaSIS$({\alpha=.75})$ &4& $60.5~ (176.3)$ & 0.87 \NN
% &NIS&4& $742.0~ (365.5)$ & 0.23 \NN
% &SIRS &4& $495.0~( 500.0)$ & 0.77 \ML
%%%%%%%%%%%%%%%%%%%%%%%%%%%%
&(1b) &QaSIS $({\alpha=0.50})$ & \phantom{0}4&4 (0) & 1.00 \\
& &QaSIS $({\alpha=0.75})$& \phantom{0}4& 4 (0) & 1.00 \\
&&NIS& \phantom{0}4& 4 (0) & 1.00 \\
&&SIRS &\phantom{0}4& 6 (9) & 0.99 \\[3pt]
%%%%%%%%%%%%%%%%%%%%%%%%%%%%%
&(1c) &QaSIS$({\alpha=0.50})$& \phantom{0}4 & 4 (0) & 1.00\\
& &QaSIS$({\alpha=0.75})$& \phantom{0}4 & 4 (0)& 1.00\\
&&NIS &\phantom{0}4 & 6 (79) &0.83 \\
&&SIRS &\phantom{0}4 & 7 (14) &0.98\\[6pt]
%__________________________________________%
Example~\ref{ex2}& &QaSIS $({\alpha=0.50})$&\phantom{0}5 & $6~(2)$ & 1.00 \\
&&QaSIS $({\alpha=0.75})$ &\phantom{0}8& 18 (24) & 0.96 \\
%%&&NIS(L)& 5&$1592~(615)$ & 0.04 \NN
&&NIS &\phantom{0}8& 1726 (511) & 0.22 \\
&&SIRS &\phantom{0}8 & 18 (16) & 0.97 \\[6pt]
%%\cmidrule(r){2-6}
%&2b&QaSIS$({\alpha=.50})$&5 & $7~(3)$ & 1.00 \NN
%&&QaSIS$({\alpha=.75})$ &8& $56~ (72)$ & 0.79 \NN
%&&NIS(L)& 5&$1672~ (522)$ & 0.02 \NN
%&&NIS(LS)&8& $1747~(431)$ & 0.15 \NN
%&&SIRS &8 & $12~(5)$ & 0.99 \ML
%__________________________________________%
Example~\ref{ex3}&(3a) &QaSIS $({\alpha=0.50})$ &\phantom{0}2& 3 (2) & 1.00\\
& &QaSIS $({\alpha=0.75})$&15& 153 (207) & 0.89 \\
% &NIS(L) &2& $1307.5~( 900.5)$ & 0.05\NN
& &NIS&15& 3117 ( 4071) & 0.50\\
&&SIRS &15& 698 (1140) & 0.89\\[3pt]
%%%%%%%%%%%%%%%%%%%%%%%%%%%%%%%%%%
&(3b) & QaSIS $({\alpha=0.50})$ &\phantom{0}2& 2 (1) & 1.00\\
& &QaSIS $({\alpha=0.75})$ &15& 88 (542) & 0.92 \\
%&NIS(L)&15& $1260.50~ (902.50)$ & 0.06 \NN
&&NIS& 15&4166 (1173) & 0.25 \\
&&SIRS &15& 29 (21) & 1.00 \\
\hline
\end{tabular*}
\tabnotetext[]{}{$p^{*}$: the number of truly active variables.}
\end{table}

The median value of $\mathcal{R}$ (with IRQ in the parenthesis) and
the average value of
$\mathcal{S}$ for QaSIS, NIS and SIRS are summarized in Table~\ref{tab1}. For
QaSIS, we report results for two quantiles $\alpha=0.5$ and 0.75.
We observe the following from Table~\ref{tab1}: (i)~The $L_2$ norm based NIS procedure
exhibits the best performance when the random error has a normal
distribution, but its performance deteriorates
substantially for heavy-tailed or heteroscedasitic errors (Examples
\ref{ex1}--\ref{ex2}). (ii)
We observe that in case (1a) where $\rho=0$, no
method works really well in terms of the minimum model size. This is because
the independent signals work against the marginal effect estimation as
accumulated noise, thus masking the relatively weak signals from $X_3$
and $X_4$ in this model.
%behave as accumulated noise when estimating the marginal effects};
%SIRS performs poorly when the covariates
%do not follow the multivariate normal distribution (\bfblue{we will need
%to remove this statement as we have removed the original example 1a,
%where the covariate
%distribution is uniform. Maybe we should add it back?}), but is robust
%to heavy-tailed and heteroscedastic %errors, and
(iii) In Example~\ref{ex3} where the model has a more complex structure,
QaSIS is effective in identifying the number of active variables at
different quantiles;
while the performance of SIRS depends on the functional form.
Overall, our simulations for the complete data case demonstrate that
the performance of QaSIS is on par with or better than that of NIS and
SIRS for a variety of distributions of covariates and errors.
%______________________________________________________%

Variable screening with censored responses has received little
attention in the literature.
In Example~\ref{ex4} below, we compare the quantile-adaptive nonparametric
marginal screening
procedure proposed in Section~\ref{sec4} with the Cox model based marginal
screening procedure
[Cox(SIS)] of Fan, Feng and Wu (\citeyear{F12}) and a naive procedure treating
the censored data as complete
and then applying to the QaSIS procedure (denoted by Naive).

\begin{example}[(Censored responses)]\label{ex4}
We consider a case in which the
latent response variable $Y_i$
is generated using the same setup as in case (1b). Let $Y_i^*=\min
(Y_i,C_i)$, where
the censoring time $C_i$ is generated from a 3-component normal mixture
distribution $0.4N(-5,4)+0.1N(5,1)+\break0.5N(55,1)$.
The censoring probability is about 45\%. Due to the high censoring
rate, the performance of the variable screening procedures is
investigated at the median and the 0.25 quantile.
\end{example}

%We refer to our modified QaSIS($\alpha$) for the censored data as
%QaSIS (Surv;$\alpha$).
%Figure~\ref{figure} compares QaSIS(Surv;.5) and naive estimator where
%the censoring was completely ignored.
%Thus, $Y_i$ is treated as $Y^*_i$ and fitted using QaSIS.
%The short dashed line is the fitted regression for the latent $Y_i$.
%The long dashed and solid lines are for QaSIS (Surv;.5) and Naive,
%respectively.
%Censored observations are black dotted.
%Figure~\ref{figure} demonstrates that when ignoring censoring, the
%bias was quite visible, but
%QaSIS(Surv) can effectively handle the censoring.

Table~\ref{tab2} summarizes the simulations results based on 100 runs.
QaSIS substantially outperforms both Naive and Cox(SIS).
Under-performance of Cox(SIS) can be attributed to the fact that the
proportional hazards assumption is not satisfied in this example.
Table~\ref{tab2} also includes the LQaSIS procedure which will be discussed in
Section~\ref{sec7.1}.
%%%%%%%%%%%%%%%%%%%%%%%%%%%%%%%%%%%%%%%%%%%%%%%%%%%%%

%
\begin{table}
\caption{Simulation results for Example \protect\ref{ex4} ($n=400$, $p=1000$).
The numbers reported are the median of $\mathcal{R}$ [with
interquartile range (IQR) given in parentheses] and $\mathcal{S}$}\label{tab2}
\begin{tabular*}{\textwidth}{@{\extracolsep{\fill}}lccc@{}}
\hline
\textbf{Method}&
\multicolumn{1}{c}{$\bolds{p^*}$}&
\multicolumn{1}{c}{$\bolds{\mathcal{R}}$\textbf{(IQR)}}
&\multicolumn{1}{c@{}}{{$\bolds{\mathcal{S}}$}}\\
\hline
QaSIS ($\alpha=0.50$)& 4&4 (2) & 0.99\\
QaSIS ($\alpha=0.25$)& 4&5 (22) & 0.96\\
LQaSIS ($\alpha=0.50$)& 4&4 (6) & 0.98\\
LQaSIS ($\alpha=0.25$)& 4& 5 (13) & 0.98\\
Naive ($\alpha=0.50$)& 4&254 (497) & 0.74\\
Naive ($\alpha=0.25$)& 4&792 (351) & 0.14\\
Cox (SIS) &4& 190 (655) & 0.70\\
\hline
\end{tabular*}
\tabnotetext[]{}{$p^{*}$: the number of truly active variables.}
\end{table}
%

% \ctable[ caption={
% Simulation results for Example 4.
%The numbers reported are the median of $\mathcal{R}$ with
%interquartile range (IQR) given in parentheses, and %$\mathcal{S}$.
%},
%label=t2, pos=h!, ]{cllll} { \tnote[]{Note: Same caption as in Table~
%(IQR) &$\mathbf{\mathcal{S}}$\ML
% &QaSIS($\alpha=.75$) &4& $46.0~(114.3)$ & 0.90\NN
% &Naive($\alpha=.50$)& 4&$155.0~ (313.8)$ & 0.79 \NN
% &Naive($\alpha=.75$)& 4&$65.0~ (178.0)$ & 0.87 \NN
% &Cox(SIS) &4& $527.0~(521.25)$ & 0.75 \LL
% }

%Figure~\ref{figure} compares QaSIS(Surv;.5) and naive estimator where
%the censoring was completely ignored.
%(\bfblue{I temporarily removed the figure in an earlier version of the
%draft as it concerns the estimation performance.})

%method.
% The short and long dashed lines are the fitted $Y^*_i$ using QaSIS
%and QaSIS (Surv), respectively.
% The solid line is for the Naive method.
%Censored observations are black dotted.}\label{figure}

%%%%%%%%%%%%%%%%%%%%%%%%%%%%%%%%%%%%%%%%%%%%%%%%%%%%%%%%%%%%%%%%%%%%%
%s6 #&#
\section{Real data analysis}\label{sec6}
We illustrate the proposed screening method on the diffuse large-B-cell
lymphoma (DLBCL) microarray data of \citet{Ro}. The data
set contains the survival times of 240 patients and the
gene expression measurements of 7399 genes for each patient.
%The outcome in the study was the survival time, and gene expression
%measurements of $p=7399$ features were used %as predictors.
%Since high dimensionality of this data set poses challenges in
%handling and analyzing a vast amount of %information, and the cause of
%a clinical phenotype might involve only a small set of genes, it is
%desirable to %implement dimension reduction methods prior to fitting
%the data into the proposed model.
The gene expression measurements
for each gene are standardized to have mean zero and variance one.
To assess the predictive performance of the proposed method, we divide
the data set into a training set with $n_1=160$ patients and a testing
set with remaining $n_2=80$ patients, in the same way as \citet{Bair} did.
The index of the training set is available from
\url{http://www-stat.stanford.edu/\textasciitilde tibs/superpc/staudt.html}.

Nearly half of the survival time data are censored, so we focus
our attention on two quantile levels $\alpha=0.2$ and 0.4 that
represent the effects of gene expression on the sub-population of
patients with poor prognosis. We apply the proposed QaSIS method to the
training data to
select $[n_1/\log(n_1)]=31$ genes, which is followed up by
a variable selection procedure based on additive quantile regression
with the SCAD penalty
[\citet{FL}] to find two top genes. As with the empirical studies
in the simulation study, we use three internal knots of $N=3$. Because
almost all the censoring occurs above the 0.4th quantile, we do not
need to re-weight the censored observations for the low quantiles we
are considering in this example. Based on the two selected genes at
each $\alpha$, we estimate the corresponding quantile function. The
estimated quantile function from the training set is then used to
calculate risk scores for each patient in the testing data set. If
$Y_i$ is the survival time of the $i$th patient in the training set,
with $s_i$ as the predicted risk score, we expect the $\alpha$th
quantile of $Y_i$ given $s_i$ to have a significant (and positive) slope.

For the purpose of comparison, we also follow the same analysis path
but replace QaSIS by the sure independence screening for Cox models,
SIS(Cox), of \citet{F12} and the SCAD-penalized Cox
regression to select two genes. The risk scores are then calculated
based on the linear index for the Cox model. Table~\ref{real:table3}
summarizes the slope coefficients of regressing survival times on risk
scores in the training set based on the censored quantile regression of
\citet{Portnoy-2003}. It is clear that the analysis based on QaSIS has the
desired predictive power, where the 0.2 and 0.4 quantiles of
survival time for the testing data set are significantly associated
with the predicted risk scores, but the analysis based on SIS(Cox) did
not make it. If we regress the survival time on the risk scores on the
training data, we would get coefficients of exactly 1.0 under QaSIS,
but it would not have validation power.

%
%t1 #&#
\begin{table}
\tablewidth=250pt
\caption{Estimated slope coefficients (and $p$-values) for survival
time versus the risk score at $\alpha$th quantile obtained by QaSIS
and SIS (Cox)}
\label{real:table3}
\begin{tabular*}{250pt}{@{\extracolsep{\fill}}lcc@{}}
\hline
\textbf{Method}& \textbf{Estimated coefficients} &$\bolds{p}$\textbf{-value}\\
\hline
QaSIS ($\alpha=0.4$)&0.93 &\textbf{0.02}\\
QaSIS ($\alpha=0.2$)&0.53&\textbf{0.04}\\[3pt]
%----------------------------------------------------------------------%
SIS (Cox) ($\alpha=0.4$)&0.16 & 0.62\\
SIS (Cox) ($\alpha=0.2)$&0.17&0.62\\
\hline
\end{tabular*}
\end{table}

We examined the two genes selected by QaSIS. At $\alpha=0.4$, their
GeneIDs are 31981 (known as AA262133, septin 1) and 17902 (AA284323,
glutathione synthetase). Both genes belong to the known Proliferation
signature group in the study of DLBCL by \citet{Ro}. Gene
AA262133 was also ranked very high by \citet{Li} using the
partial likelihood-based scores. At $\alpha=0.2$, the two selected
genes have IDs 31585 (known as NM 018518, MCM10 minichromosome
maintenance deficient 10) and 33014 (AA769543, Hypothetical protein
MGC4189). We find that Gene 31858 also belongs to the known
Proliferation signature group. Gene 33014 was identified as an
interesting candidate by \citet{Li}.
We did not find any signature group associated with it, but it seems
quite related to the lower tail of the survival time distribution, and
is worth further investigation.

%%%%%%%%%%%%%%%%%%%%%%%%%%%%%%%%%%%%%%%%%%%%%%%%%%%%%%%%%%%%

%s7 #&#
\section{Discussions}\label{sec7}
%s7.1 #&#
\subsection{Further extension on screening with censored responses}\label{sec7.1}
The assumption that the censoring distribution
does not depend on the covariates is popular in regression analysis of
survival data.
%In particular, it is often adopted in the so-called administrative
%censoring case, where the entry time of a study
%is staggered and the censoring occurs if a patient is still alive
%when the study ends at a prespecified time $T$.
%Despite that the independence assumption is common and
%the proposed procedure in Section 4 displays a certain degree of
%robustness
%(see additional simulation results in the supplementary file),
It can be further relaxed. Assume that $Y_i$ and $C_i$ are
conditionally independent
given $X_{ij}$, $j=1,\ldots,p$.
Let $G(t|X_{ij})=P(C_i>t|X_{ij})$ be the conditional survival function
of $C_i$ given $X_{ij}$.
Let $\widehat{G}(t|X_{ij})$ be the local Kaplan--Meier estimator of
$G(t|X_{ij})$
[e.g., \citet{Beran-1981} and Gonzalez-Manteiga and Cadarso-Suarez (\citeyear{GMCS-1994})].
More specifically,
%the local Kaplan--Meier estimator is obtained based on $\{Y_i^*, X_{ij},
%
%e7.1 #&#
\begin{eqnarray}
\label{G2hat} \widehat{G}(t|x) &=& \prod_{j=1}^n
\biggl\{1 - \frac{B_{nj}(x)}{\sum_{k=1}^n I(Y_k^*\geq Y_j^*)B_{nk}(x)} \biggr\}^{I(Y_j^*
\leq{t}, \delta_j =0)},
\end{eqnarray}
where $
B_{nk}(x)
=K(\frac{x-{x_k}}{h_n})/\{\sum_{i=1}^nK(\frac{x-{x_i}}{h_n})\},
k=1,\ldots,n,
$
are the Nadaraya--Watson weights,
$h_n$ is the bandwidth and $K(\cdot)$ is a density function.
We consider estimating $\vbeta_{0j}$ using the
locally weighted censored quantile regression, that is,
\[
\widetilde{\vbeta}_j^{c}=\mathop{\argmin}_{\vbeta\in\mathbb{R}^N}\sum
_{i=1}^n\frac{\delta_i}{\widehat{G}(Y_i^*|X_{ij})}
\rho_{\alpha
}\bigl(Y_i^{*}-\bolds{
\pi}(X_{ij})^T\vbeta\bigr).
\]
Let
$
\widetilde{f}_{nj}^{c}(t)=\bolds{\pi}(t)^T\widetilde
{\vbeta}_{j}^{c}-F_{\mathrm{KM},n}^{-1}(\alpha)
$
and define $
\widetilde{M}_{\alpha}^{c}=\{1\leq j\leq p\dvtx \Vert \widetilde
{f}_{nj}^{c}\Vert_n^2\geq\eta_n^{c}\}
$
where $\eta_n^{c}$ is a predefined threshold value. We refer to this
new procedure
as LQaSIS, whose numerical performance is reported in Table~\ref{tab2} and two other
examples in the supplemental material [He, Wang and Hong (\citeyear{HWH})].

We assume, instead of (C6):
\begin{longlist}[(C6$^{\prime}$)]
\item[(C6$^{\prime}$)] $\inf_x P(t \leq Y_i \leq C_i|x) \geq \tau_0 > 0$
for some positive constant $\tau_0$ and any $t \in [0, T]$, where $T$ denotes
the maximum follow-up time. $G(t|x)$ has first derivatives with respect to $t$,
which is uniformly bounded away from infinity; and
$G(t|x)$ has bounded (uniformly in $t$) second-order partial
derivatives with respect to $x$. Furthermore,
$t_0\leq\sup\{t\dvtx G(t|x)>0\}\leq t_1$ uniformly in $x$ for some
positive constants $t_0$ and $t_1$,
and $\sup\{t\dvtx P(Y>t|x)>0\}\geq\sup\{t\dvtx G(t|x)>0\}$ almost surely for
$x$.
\end{longlist}

Then Theorem~\ref{main2} can be extended as follows.
%
%th7.1 #&#
\begin{theorem}\label{main2_local}
Assume conditions \textup{(C1)--(C5)}, \textup{(C6$^{\prime}$)} and \textup{(C7)} are satisfied,
$\tau<1/4$,
$nh^3\ra\infty$, $N^{3}n^{2\tau-1}=o(1)$,
%$n^{2\tau-3/4}(\log n)^{3/4}h^{-3/4}=o(1)$,
$N^2n^{2\tau-1}(\log n)^2h^{-1}=o(1)$
and $(N+n^{\tau})n^{\tau}h^2=o(1)$.
If we take $\nu_n=\delta^* n^{-\tau}$ with $\delta^*\leq c_1/16$, then:
%then $n^{2\tau-3/4}(\log n)^{3/4}h^{-3/4}=o(1)$).}
\begin{longlist}[(1)]
\item[(1)] there exist positive constants $\delta_3$ and $\delta_4$ such that
\begin{eqnarray*}
P \bigl(M_{\alpha}\subset\widehat{M}_{\alpha}^{c} \bigr)
\geq1-S_{\alpha} \bigl\{17\exp \bigl(-\delta_3n^{1-4\tau}
\bigr) +12N^2\exp \bigl(-\delta_{4}N^{-3}n^{1-2\tau}
\bigr) \bigr\},
\end{eqnarray*}
for all $n$ sufficiently large. Especially,
$
P (M_{\alpha}\subset\widehat{M}_{\alpha}^{c} )\ra1
$ as $n\ra\infty$.

\item[(2)]
\begin{eqnarray*}
&& P \bigl(\bigl|\widehat{M}_{\alpha}^{c}\bigr|\leq2N^2n^{\tau}
\lambdamax (\vSigma )/\delta^*  \bigr)
\\
&&\qquad \geq 1-p \bigl\{17\exp \bigl(-b_7n^{1-4\tau} \bigr)
+12N^2\exp \bigl(-b_{8}N^{-3}n^{1-2\tau}
\bigr) \bigr\},
\end{eqnarray*}
for all $n$ sufficiently large. Especially,
$P (|\widehat{M}_{\alpha}^{c}|\leq2N^2n^{\tau}\lambdamax(\vSigma )
/\delta^*   )\ra1$
as \mbox{$n\ra\infty$}.
\end{longlist}
\end{theorem}

The proof of Theorem~\ref{main2_local} is given in the supplemental
material [He, Wang and Hong (\citeyear{HWH})].

%s7.2 #&#
\subsection{Limitations and other issues}\label{sec7.2}

We have not investigated the problem of adaptively selecting the number
of basis functions
in this paper for two reasons: (1) although adaptive tuning is
possible, it will significantly increase the computational time
as it needs to be done for each covariate separately; (2) optimal
estimation is not the goal for marginal
screening, instead consistent estimation generally suffices.
Empirically, we find that 3 or 4 internal knots are generally enough to
flexibly approximate many smooth functions typically seen in practice.

Marginally unimportant but jointly important variables may not be
preserved in marginal screening.
This is a well-recognized weakness of all existing marginal screening
procedures.
Iterative procedures may be helpful to a certain degree for this
problem [\citet{F2}].
In the same spirit, we find that in practice
a slightly modified QaSIS helps in situations where a dominating
variable increases the error variance of the marginal regression model
for other variables and hence mask the significance of other variables.
If the top ranked variable is dominating, then the modified QaSIS
removes its effects from $\vY$ first and screen the remaining
variables again.

The way we define the set of active variables can be considered
as a nonparametric approach in the sense that we consider directly the
conditional
quantile function without a specific model structure. In real life
high-dimensional data analysis,
the knowledge needed for an
appropriate model specification is often inadequate. Using a
misspecified model
to perform variable selection is likely to produce misleading results.
We propose
to separate variable screening and model building, where a
nonparametric approach is applied
to screen high-dimensional variables and then followed by sensible
model building
in the second stage in a lower-dimensional space. We expect that
this model-free approach to variable screening to
gain momentum in ultra-high dimensional learning,
see, for example, the work of \citet{LZZ} on distance
correlation based
screening.

%s8 #&#
\section{Technical proofs}\label{sec8}
We present the proof for the random censoring case, as this is the more
challenging scenario.
The proof for the complete data case and that for Lemma~\ref{lemma1}
are given in the supplemental material
[He, Wang and Hong (\citeyear{HWH})].\vadjust{\goodbreak}

To establish the sure independence property, a key step is to obtain an
exponential tail probability
bound for
%
%e8.1 #&#
\begin{eqnarray}
\label{target} P \Bigl(\max_{1\leq j\leq p} \bigl|\bigl\Vert \widehat {f}_{nj}^c
\bigr\Vert_n^2-\Vert f_{nj}\Vert^2 \bigr|\geq
Cn^{-\tau} \Bigr),
\end{eqnarray}
where $C$ is any positive constant.
We have
\begin{eqnarray*}
\bigl\Vert \widehat{f}_{nj}^c\bigr\Vert_n^2&=&
\widehat{\vbeta}_j^{cT} \bigl(\mathbb{P}_n
\bolds{\pi}(X_j)\bolds{\pi}(X_j)^T\bigr)
\widehat{\vbeta}_j^c -2F_{\mathrm{KM},n}^{-1}(
\alpha)n^{-1}\bigl(\mathbb{P}_n\bolds{\pi
}(X_j)\bigr)^T\widehat{\vbeta}_j^{c}\\
&&{}+\bigl(F_{\mathrm{KM},n}^{-1}(\alpha)\bigr)^2,
\\
\Vert f_{nj}\Vert^2&=&\vbeta_{0j}^T
\bigl(\E\bolds{\pi}(X_j)\bolds{\pi }(X_j)^T
\bigr)\vbeta_{0j} -2Q_{\alpha}(Y) \bigl(E\bolds{
\pi}(X_j)\bigr)^T\vbeta_{0j}+
\bigl(Q_{\alpha}(Y)\bigr)^2,
\end{eqnarray*}
where $\mathbb{P}_n\bolds{\pi}(X_j)\bolds{\pi}(X_j)^T=n^{-1} \sum_{i=1}^n\bolds{\pi}(X_{ij})\bolds{\pi}(X_{ij})^T$,
$\mathbb{P}_n\bolds{\pi}(X_j)=n^{-1}\times\break \sum_{i=1}^n\bolds{\pi
}(X_{ij})$, and
$\E\bolds{\pi}(X_j)\bolds{\pi}(X_j)^T$
denotes the expectation of $\bolds{\pi}(X_j)\bolds{\pi}(X_j)^T$
under the true distribution of $X_j$.
Note that
\begin{eqnarray*}
&&\bigl\Vert \widehat{f}_{nj}^c\bigr\Vert_n^2-
\Vert f_{nj}\Vert^2
\\
&&\qquad= \bigl(\widehat{\vbeta}_j^c-\vbeta_{0j}
\bigr)^T\bigl(\mathbb{P}_n\bolds{\pi }(X_j)
\bolds{\pi}(X_j)^T\bigr) \bigl(\widehat {
\vbeta}_j^c-\vbeta_{0j}\bigr)
\\
&&\quad\qquad{}+2\bigl(\widehat{\vbeta}_j^c-\vbeta_{0j}
\bigr)^T\bigl(\mathbb{P}_n\bolds{\pi }(X_j)
\bolds{\pi}(X_j)^T\bigr)\vbeta_{0j}
\\
&&\quad\qquad{}+\widehat{\vbeta}_j^c \bigl(\mathbb{P}_n
\bolds{\pi}(X_j)\bolds {\pi}(X_j)^T-\E
\bolds{\pi}(X_j)\bolds{\pi}(X_j)^T \bigr)
\vbeta_{0j}
\\
&& \quad\qquad{}-2F_{\mathrm{KM},n}^{-1}(\alpha) \bigl[\mathbb{P}_n
\bolds{\pi }(X_j)^T\widehat{\vbeta}_j^{c}-
E\bolds{\pi}(X_j)^T\vbeta_{0j} \bigr]
\\
&& \quad\qquad {}+2 \bigl[Q_{\alpha}(Y)-F_{\mathrm{KM},n}^{-1}(\alpha) \bigr]
\bigl(E\bolds{\pi}(X_j)^T\vbeta_{0j} \bigr)\\
&&\qquad\quad{} +
\bigl[\bigl(F_{\mathrm{KM},n}^{-1}(\alpha)\bigr)^2-
\bigl(Q_{\alpha}(Y)\bigr)^2 \bigr]
\\
&&\qquad\stackrel{\Delta} {=}  \sum_{k=1}^6S_{jk},
\end{eqnarray*}
where the definition of $S_{jk}$ is clear from the context.
From the argument of Lemma~\ref{lemma1}, $E(\bolds{\pi
}(X_j)^T\vbeta_{0j})$ is uniformly bounded in $X_j$
and by Lemma~\ref{lemma8}(4) below, we have $|S_{j5}|=O(n^{-1/2}(\log
n)^{1/2})=o(n^{-\tau})$ almost surely.
Similarly, $|S_{j6}|=O(n^{-1/2}(\log n)^{1/2})=o(n^{-\tau})$ almost
surely. Therefore, for all $n$
sufficiently large,
\begin{eqnarray*}
&&P \Bigl(\max_{1\leq j\leq p} \bigl|\bigl\Vert \widehat {f}_{nj}^c
\bigr\Vert_n^2-\Vert f_{nj}\Vert^2 \bigr|\geq
Cn^{-\tau} \Bigr)
\\
&&\qquad\leq P \Biggl(\max_{1\leq j\leq p}\sum_{k=1}^4
|S_{jk} |\geq Cn^{-\tau}/2 \Biggr)
\\
&&\qquad \leq \sum_{k=1}^4P \Bigl(
\max_{1\leq j\leq p} |S_{jk} |\geq Cn^{-\tau}/8 \Bigr).
\end{eqnarray*}
In the following, we shall provide details on deriving exponential tail
bound for
$P (\max_{1\leq j\leq p} |S_{jk} |\geq Cn^{-\tau}/8 )$.

%s8.1 #&#
\subsection{Properties of the spline basis}\label{sec8.1}

First, we recall some useful properties of the basis vector $\bolds
{\pi}(t)=(B_1(t),\ldots,B_N(t))^T$.
Zhou, Shen and Wolfe (\citeyear{Zhou}) established that there exist two positive constants
$b_1$ and $b_2$ such that
%
%e8.2 #&#
\begin{eqnarray}
\label{eigenbound} b_1N^{-1}&\leq&\lambdamin\bigl(\E
\bolds{\pi}(X_j)\bolds{\pi}(X_j)^T\bigr)
\leq\lambdamax\bigl(\E\bolds{\pi}(X_j)\bolds{\pi}(X_j)^T
\bigr)
\nonumber
\\[-8pt]
\\[-8pt]
\nonumber
&\leq& b_2N^{-1}\qquad  \forall j,
\end{eqnarray}
where $\lambdamin$ and $\lambdamax$ denote the smallest eigenvalue
and the largest eigenvalue, respectively.

\citet{Stone} established that there exists a positive constant $b_3$
such that
%
%e8.3 #&#
\begin{eqnarray}
\label{bbound4} E\bigl(B_k^2(X_{ij})\bigr)\leq
b_3N^{-1},\qquad  1\leq k\leq N, 1\leq i\leq n, 1\leq j\leq p.
\end{eqnarray}
Similar result can be found in \citet{He2}.

%le8.1 #&#
\begin{lemma}\label{lemma3} Let $\mathbb{P}_n\bolds{\pi}(X_j)\bolds
{\pi}(X_j)^T=n^{-1}\sum_{i=1}^n\bolds{\pi}(X_{ij})\bolds{\pi
}(X_{ij})^T$ and
$\vD_j=\mathbb{P}_n\bolds{\pi}(X_j)\bolds{\pi}(X_j)^T-\E\bolds
{\pi}(X_j)\bolds{\pi}(X_j)^T$.
\begin{longlist}[(1)]
\item[(1)] There exists a positive constant $c_4$ such that for all $n$
sufficiently large
\begin{eqnarray}
\label{bbound} \qquad P \bigl(\lambdamax\bigl(\mathbb{P}_n\bolds{
\pi}(X_j)\bolds{\pi}(X_j)^T\bigr)
\geq(b_2+1)N^{-1} \bigr) \leq2N^2\exp
\bigl(-c_4nN^{-3} \bigr),
\end{eqnarray}
\item[(2)] For any $c_5>0$, there exists a positive constant $c_6$ such that
for all $n$ sufficiently large
\begin{eqnarray}
\label{bbound2} &&
P \bigl(\max \bigl(\bigl|\lambdamax(\vD_j)\bigr|, \bigl|\lambdamin(
\vD_j)\bigr| \bigr) \geq c_5N^{-1}n^{-\tau}
\bigr)
\nonumber
\\[-8pt]
\\[-8pt]
\nonumber
&&\qquad \leq2N^2\exp \bigl(-c_6N^{-3}n^{1-2\tau}
\bigr).
\end{eqnarray}
\end{longlist}
\end{lemma}

\begin{pf} The proof is an extension of that for Lemma 5 in \citet{F1} which proved similar results for the smallest eigenvalue.
First, for any two symmetric
matrices $\vA$ and $\vB$, we have $\lambdamax(\vA+\vB)\leq
\lambdamax(\vA)+\lambdamax(\vB)$.
This implies that $\lambdamax(\vA)-\lambdamax(\vB)\leq\lambdamax
(\vA-\vB)$
and $\lambdamax(\vB)-\lambdamax(\vA)\leq\lambdamax(\vB-\vA)$. Thus
\begin{eqnarray*}
\bigl|\lambdamax(\vA)-\lambdamax(\vB) \bigr|\leq\max \bigl\{ \bigl|\lambdamax(\vA-\vB)\bigr|, \bigl|
\lambdamax(\vB-\vA)\bigr| \bigr\}.
\end{eqnarray*}
Applying the above inequality, we have
%
%e8.4 #&#
\begin{eqnarray}
\label{cat1} && \bigl|\lambdamax\bigl(\mathbb{P}_n\bolds{
\pi}(X_j)\bolds{\pi }(X_j)^T\bigr)-
\lambdamax\bigl(\E\bolds{\pi}(X_j)\bolds{\pi}(X_j)^T
\bigr) \bigr|
\nonumber
\\[-8pt]
\\[-8pt]
\nonumber
&&\qquad \leq\max \bigl\{\bigl|\lambdamax(\vD_j)\bigr|, \bigl|\lambdamax(-
\vD_j)\bigr| \bigr\}.
\end{eqnarray}
For any $N$-dimensional vector
$\mathbf{a}=(a_1,\ldots,a_N)^T$ satisfying $\Vert \mathbf{a}\Vert =1$, we have
$
|\mathbf{a}^T\vD_j\mathbf{a}|\leq\Vert \vD_j\Vert_{\infty} (\sum_{i=1}^N|a_i| )^2\leq N\Vert \vD_j\Vert_{\infty},
$
where $\Vert \vD_j\Vert_{\infty}$ is the sup norm of the matrix $\vD_j$.
%=\max_{1\leq i,l\leq N}|\vD_j^{(i,l)}|$
%with $\vD_j^{(i,l)}$ being the $(i,l)$th entry of $\vD_j$.
Thus $\lambdamax(\vD_j)=\max_{\Vert \mathbf{a}\Vert =1}\mathbf{a}^T\vD_j\mathbf{a}\leq N\Vert \vD_j\Vert_{\infty}$.
Also\break $\lambdamax(\vD_j)=-\min_{\Vert \mathbf{a}\Vert =1} (-\mathbf
{a}^T\vD_j\mathbf{a} )\geq-N\Vert \vD_j\Vert_{\infty}$.\vadjust{\goodbreak}
Thus $|\lambdamax(\vD_j)|\leq N\Vert \vD_j\Vert_{\infty}$. Similarly, we have
$|\lambdamax(-\vD_j)|\leq N\Vert \vD_j\Vert_{\infty}$.
Following (\ref{cat1}) and using the result on the smallest eigenvalue
of $\vD_j$ [Fan, Feng and Song (\citeyear{F1})],
we have
%
%e8.5 #&#
\begin{eqnarray}
\label{bbound5} &&\bigl |\lambdamax\bigl(\mathbb{P}_n\bolds{
\pi}(X_j)\bolds{\pi }(X_j)^T\bigr)-
\lambdamax\bigl(\E\bolds{\pi}(X_j)\bolds{\pi}(X_j)^T
\bigr) \bigr|
\nonumber
\\[-8pt]
\\[-8pt]
\nonumber
&&\qquad\leq\max\bigl(\bigl|\lambdamax(\vD_j)\bigr|, \bigl|\lambdamin(
\vD_j)\bigr|\bigr)\leq N\Vert \vD_j\Vert_{\infty}.
\end{eqnarray}

As in \citet{F1}, applying Bernstein's inequality to
each entry of
$\vD_j$, it can be shown that $\forall\delta>0$,
%
%e8.6 #&#
\begin{equation}
\label{bbound7} P \bigl(N\Vert \vD_j\Vert_{\infty}\geq N\delta/n \bigr)
\leq2N^2\exp \biggl\{-\frac{\delta^2}{2(b_3nN^{-1}+\delta/3)} \biggr\}.
\end{equation}

To prove (\ref{bbound}), we use the bound in~(\ref{eigenbound}),
apply the inequality in~(\ref{bbound5}) and
take $\delta=N^{-2}n$ in~(\ref{bbound7}). This gives
\begin{eqnarray*}
P \bigl(\lambdamax\bigl(\mathbb{P}_n\bolds{\pi}(X_j)
\bolds{\pi}(X_j)^T\bigr) \geq(b_2+1)N^{-1}
\bigr) \leq2N^2\exp \bigl(-c_4N^{-3}n \bigr),
\end{eqnarray*}
for some positive constant $c_4$ for all $n$ sufficiently large.

To prove (\ref{bbound2}), we apply the inequality in (\ref{bbound5}) and
take $\delta=c_5N^{-2}n^{1-\tau}$ in~(\ref{bbound7}). This gives
\begin{eqnarray*}
P \bigl(\max \bigl(\bigl|\lambdamax(\vD_j)\bigr|, \bigl|\lambdamin(
\vD_j)\bigr| \bigr) \geq c_5N^{-1}n^{-\tau}
\bigr) \leq2N^2\exp \bigl(-c_6N^{-3}n^{1-2\tau}
\bigr),
\end{eqnarray*}
for some positive constant $c_6$ for all $n$ sufficiently large.
\end{pf}

\subsection{\texorpdfstring{An exponential tail probability bound for $\Vert {\widehat{\vbeta}}_j^c-\vbeta_{0j}\Vert $}
{An exponential tail probability bound for ||beta jc - beta 0j ||}}\label{sec8.2}

Let
\begin{eqnarray*}
B_n(\vbeta)&=&n^{-1}\sum_{i=1}^n
\delta_i\bigl[\widehat{G}\bigl(Y_i^*\bigr)
\bigr]^{-1} \bigl[\rho_{\alpha}\bigl(Y_i^*-\bolds{
\pi}(X_{ij})^T\vbeta \bigr)-{\rho_{\alpha}
\bigl(Y_i^*\bigr)} \bigr],
\\
B(\vbeta)&=&\E \bigl\{\delta_i\bigl[G\bigl(Y_i^*\bigr)
\bigr]^{-1} \bigl[\rho_{\alpha
}\bigl(Y_i^*-\bolds{
\pi}(X_{ij})^T\vbeta\bigr)-{\rho_{\alpha
}
\bigl(Y_i^*\bigr)} \bigr] \bigr\}.
\end{eqnarray*}
Then $\widehat{\vbeta}_{j}^c=\argmin_{\vbeta\in\mathbb
{R}^N}B_n(\vbeta)$.
Applying the iterative expectation formula, we have
\begin{eqnarray*}
B(\vbeta)&=&\E \bigl\{\E \bigl\{I(Y_i\leq C_i)
\bigl[G(Y_i)\bigr]^{-1} \bigl[\rho_{\alpha}
\bigl(Y_i-\bolds{\pi}(X_{ij})^T\vbeta\bigr) -{
\rho_{\alpha}(Y_i)} \bigr] |Y_i,X_{ij}
\bigr\} \bigr\}
\\
&=& \E \bigl[\rho_{\alpha}\bigl(Y_i-\bolds{
\pi}(X_{ij})^T\vbeta\bigr)-{\rho_{\alpha}(Y_i)}
\bigr].
\end{eqnarray*}
Hence, $\vbeta_{0j}=\argmin_{\vbeta\in\mathbb{R}^N}B(\vbeta)$.

We can bound the difference $\Vert \widehat{\vbeta}_j^c-\vbeta_{0j}\Vert $
by the difference of their respective objective functions.

%le8.2 #&#
\begin{lemma}\label{lemma5}
For any $\delta>0$,
%
%e8.7 #&#
\begin{eqnarray}
\label{bear1}\qquad &&P \bigl(\bigl\Vert \widehat{\vbeta}_j^c-
\vbeta_{0j}\bigr\Vert \geq\delta \bigr)
\nonumber
\\[-8pt]
\\[-8pt]
\nonumber
&&\qquad\leq P \biggl(\sup_{\Vert \vbeta-\vbeta_{0j}\Vert \leq\delta}
\bigl|B_n(\vbeta )-B(\vbeta) \bigr|\geq  \frac{1}{2}\inf_{\Vert \vbeta-\vbeta_{0j}\Vert =\delta} \bigl(B(\vbeta )-B(
\vbeta_{0j})\bigr) \biggr).
\end{eqnarray}
\end{lemma}

\begin{pf} This is a direct application of Lemma 2 of \citet{Hjort} making use
of the convexity of the objective function.
\end{pf}

The lower bound of the right-hand side of (\ref{bear1}) can be
explicitly evaluated for any given $\delta>0$.
This is summarized in the following lemma.

%le8.3 #&#
\begin{lemma}\label{lemma4} Let $C>0$ be an arbitrary constant.
Assume that $N^{-d}n^{\tau}=o(1)$, then
there exists a positive constant $b_4$ such that
\[
\inf_{\Vert \vbeta-\vbeta_{0j}\Vert =CN^{1/2}n^{-\tau}} \bigl(B(\vbeta )-B(\vbeta_{0j})\bigr)\geq
b_4n^{-2\tau}
\]
for all $n$ sufficiently large.
\end{lemma}

\begin{pf} We consider $\vbeta=\vbeta_{0j}+CN^{1/2}n^{-\tau}\vu$,
where $\vu\in\mathbb{R}^N$ satisfying \mbox{$\Vert \vu\Vert =1$}.
Using the identity by Knight [(\citeyear{Knight}), page 758], we have
\begin{eqnarray*}
&&B(\vbeta)-B(\vbeta_{0j})
\\
&&\qquad=\E \bigl\{\rho_{\alpha}\bigl(Y-\bolds{\pi}(X_{j})^T
\vbeta_{0j}- CN^{1/2}n^{-\tau}\bolds{
\pi}(X_{j})^T\vu\bigr)-\rho_{\alpha
}\bigl(Y-\bolds{
\pi}(X_{j})^T\vbeta_{0j}\bigr) \bigr\}
\\
&&\qquad=CN^{1/2}n^{-\tau}\E \bigl\{ \bolds{\pi}(X_{j})^T
\vu \bigl[I\bigl(Y-\bolds{\pi}(X_{j})^T
\vbeta_{0j}\leq0\bigr)-\tau \bigr] \bigr\}
\\
&&\qquad\quad{} +\E \biggl\{ \int_0^{CN^{1/2}n^{-\tau}\bolds{\pi}(X_{j})^T\vu} \bigl[I\bigl(Y-
\bolds{\pi}(X_{j})^T\vbeta_{0j}\leq s\bigr)\\
&&\hspace*{107pt}\qquad\quad{}-I\bigl(Y-\bolds{\pi}(X_{j})^T\vbeta_{0j}
\leq0\bigr) \bigr]\,ds \biggr\}
\\
&&\qquad= CN^{1/2}n^{-\tau}\E \bigl\{ \bolds{\pi}(X_{j})^T
\vu \bigl[F_{Y|X_j}\bigl(\bolds{\pi}(X_{j})^T
\vbeta_{0j}\bigr)-F_{Y|X_j}\bigl(f_j(X_j)
\bigr) \bigr] \bigr\}
\\
&&\qquad\quad{} +\E \biggl\{ \int_0^{CN^{1/2}n^{-\tau}\bolds{\pi}(X_{j})^T\vu}
\bigl[F_{Y|X_j}\bigl(\bolds{\pi}(X_{j})^T
\vbeta_{0j}+s\bigr)\\
&&\hspace*{127pt}\qquad{} -F_{Y|X_j}\bigl(\bolds{\pi}(X_{j})^T
\vbeta_{0j}\bigr) \bigr]\,ds \biggr\}
\\
&&\qquad\defined I_1+I_2.
\end{eqnarray*}
By H\"older's inequality, we have
\begin{eqnarray*}
|I_1|&\leq& CN^{1/2}n^{-\tau} \bigl(\E\bigl(\bolds{
\pi}(X_{j})^T\vu \bigr)^2 \bigr)^{1/2}
\\
&&{}\times\bigl[\E \bigl(F_{Y|X_j}\bigl(\bolds{\pi}(X_{j})^T
\vbeta_{0j}\bigr)-F_{Y|X_j}\bigl(f_j(X_j)
\bigr) \bigr)^2 \bigr]^{1/2}
\\
&\leq& CN^{1/2}n^{-\tau}O\bigl(N^{-1/2}\bigr)O
\bigl(N^{-d}\bigr)
\\
&=&O\bigl(N^{-d}n^{-\tau}\bigr),
\end{eqnarray*}
where the second inequality uses inequality (B.3) in the supplementary material
and (\ref{eigenbound}).\vadjust{\goodbreak}

Furthermore, for some $\xi$ between
$\bolds{\pi}(X_{j})^T\vbeta_{0j}+s$ and $\bolds{\pi
}(X_{j})^T\vbeta_{0j}$,
\begin{eqnarray*}
I_2&=&\E \biggl\{ \int_0^{CN^{1/2}n^{-\tau}\bolds{\pi}(X_{j})^T\vu
}f_{Y|X_j}(
\xi)s \,ds \biggr\} = O(1)\E\bigl(N^{1/2}n^{-\tau}\bolds{
\pi}(X_{j})^T\vu\bigr)^2
\\
&=&O\bigl(n^{-2\tau}\bigr)
\end{eqnarray*}
by (\ref{eigenbound}).
Note that $I_2$ is nonnegative and $I_1=o(I_2)$. Thus,
the conclusion of the lemma holds.\vspace*{-2pt}
\end{pf}

%Lemma~\ref{lemma6} and

Lemmas~\ref{lemma8}--\ref{lemma9} below provide several
useful technical results for evaluating the right-hand side of (\ref{bear1}).

%le8.4 #&#
\begin{lemma}\label{lemma8}
Assume conditions \textup{(C6)} and \textup{(C7)}. The Kaplan--Meier estimator $\widehat
{G}(t)$ satisfies:
\begin{longlist}[(1)]
\item[(1)] $\sup_{0\leq t\leq T}|\widehat{G}(t)-G(t)|=O(n^{-1/2}(\log
n)^{1/2})$ almost surely.

\item[(2)]
$
\widehat{G}(t)^{-1}-G(t)^{-1}=n^{-1}\sum_{j=1}^n\frac{\xi
(Y_j^*,\delta_j,t)}{G^{2}(t)}+R_n(t),
$
where $\xi(Y_j^*,\delta_j,t)$ are independent mean zero random
variables whose expression is given
in Theorem 1 of \citet{Lo}, and $\sup_{0\leq t\leq
T}|R_n(t)|=O(n^{-3/4}(\log n)^{3/4})$
almost surely.

\item[(3)] $\sup_{0\leq t\leq T} |\frac{1}{\widehat
{G}(t)}-\frac
{1}{G(t)} |=O(n^{-1/2}(\log n)^{1/2})$ almost surely.

\item[(4)] $\sup_{\beta_1\leq\alpha\leq\beta_2}
|F_{\mathrm{KM},n}^{-1}(\alpha)-Q_{\tau}(Y)
|=O(n^{-1/2}(\log n)^{1/2})$ almost surely.\vspace*{-2pt}
\end{longlist}
\end{lemma}
\begin{pf}
The results in (1) and (4) are
given in Lemma 3 of \citet{Lo}.
The result in (2) follows from the Taylor expansion, Theorem 1 in \citet{Lo} and the result in (1).
The proof of (3) follows Taylor expansion and~(1).\vspace*{-2pt}
\end{pf}

%le8.5 #&#
\begin{lemma}[(Massart's concentration theorem, \citeyear{Massart})]\label{lemma7}
Let $W_1,\ldots,W_n$ be independent random variables and let $\mathbb
{G}$ be a class of functions satisfying
$a_{i,g}\leq g(W_i)\leq b_{i,g}$ for some real numbers $a_{i,g}$ and
$b_{i,g}$, and for all $1\leq i\leq n$
and $g\in\mathbb{G}$.
Define
$
L^2=\sup_{g\in\mathbb{G}}\sum_{i=1}^n(b_{i,g}-a_{i,g})^2/n
$
and
$Z=\break\sup_{g\in\mathbb{G}}n^{-1} |\sum_{i=1}^n(g(W_i)-E(g(W_i))) |$.
Then for any positive $t$,
$P (Z\geq\E Z+t )\leq\exp [-\frac{nt^2}{2L^2} ]$.\vspace*{-2pt}
\end{lemma}

%le8.6 #&#
\begin{lemma}[{[Bernstein inequality for $U$-statistics,
\citet{Hoeffding}]}]\label{lemma9}
Let $U_n^2(g)$ denote the second-order $U$-statistics with kernel
function $g(t_1,t_2)$ based on the
independent random variables $Z_1,\ldots,Z_n$. Assume that the
function $g$ is bounded: $a<g<b$ for some finite constants $a$ and $b$.
If $E(g(Z_i,Z_j))=0$, $\forall i\neq j$, then $\forall t>0$,
$
P ( |U_n^2(g) |>t )\leq2 \exp (-\frac
{2kt^2}{(b-a)^2} ),
$
where $k$ denotes the integer part of $n/2$.\vspace*{-2pt}
\end{lemma}

%le8.7 #&#
\begin{lemma}\label{bear2}
%Assume that $n^{2\tau-3/4}(\log n)^{3/4}=o(1)$ and $N^2n^{2\tau-1}\log
%n=o(1)$.
Assume the conditions of Theorem~\ref{main2}.
For any $C>0$, there exist positive constants $c_7$ and $c_8$ such that
for all $n$\vadjust{\goodbreak} sufficiently\break large
\begin{eqnarray*}
P \bigl(\bigl\Vert \widehat{\vbeta}_j^c-
\vbeta_{0j}\bigr\Vert \geq CN^{1/2}n^{-\tau
} \bigr) \leq4\exp
\bigl(-c_7n^{1-4\tau} \bigr) +\exp \bigl(-c_{8}N^{-2}n^{1-2\tau}
\bigr).
\end{eqnarray*}
\end{lemma}

\begin{pf}
Following Lemmas~\ref{lemma5} and~\ref{lemma4}, there exists
some $b_4>0$ such that
for all $n$ sufficiently large,
\begin{eqnarray*}
&&P \bigl(\bigl\Vert \widehat{\vbeta}_j^c-
\vbeta_{0j}\bigr\Vert \geq CN^{1/2}n^{-\tau
} \bigr)
\\[-2pt]
&&\quad\leq P \Bigl(\sup_{\Vert \vbeta-\vbeta_{0j}\Vert \leq{CN^{1/2}n^{-\tau
}}} \bigl|B_n(\vbeta)-B(\vbeta) \bigr|\geq
b_4n^{-2\tau} \Bigr)
\\[-2pt]
&&\quad\leq P \biggl( \bigl|B_n(\vbeta_{0j})-B(\vbeta_{0j})
\bigr|\geq\frac
{1}{2}b_4n^{-2\tau} \biggr)
\\[-2pt]
&& \qquad{}+ P \biggl(\sup_{\Vert \vbeta-\vbeta_{0j}\Vert \leq{CN^{1/2}n^{-\tau}}} \bigl|B_n(\vbeta)-B_n(
\vbeta_{0j})-B(\vbeta) +B(\vbeta_{0j}) \bigr|\geq
\frac{1}{2}b_4n^{-2\tau} \biggr)
\\[-2pt]
&&\quad\defined J_1+J_2.
\end{eqnarray*}

First, we evaluate $J_1$. Let $W_i=\delta_i[G(Y_i^*)]^{-1} [\rho_{\alpha}(Y_i^*-\bolds{\pi}(X_{ij})^T\vbeta_{0j})
-\rho_{\alpha}(Y_i^*) ]$. Then
\begin{eqnarray*}
&&\hspace*{-4pt}B_n(\vbeta_{0j})-B(\vbeta_{0j})
\\[-2pt]
&&\hspace*{-4pt}\qquad=n^{-1}\sum_{i=1}^n(W_i-EW_i)
\\[-2pt]
&&\hspace*{-4pt}\qquad\quad{} +n^{-1}\sum_{i=1}^n
\delta_i \bigl[\bigl(\widehat{G}\bigl(Y_i^*\bigr)
\bigr)^{-1}-\bigl(G\bigl(Y_i^*\bigr)\bigr)^{-1}
\bigr] \bigl[\rho_{\alpha
}\bigl(Y_i^*-\bolds{
\pi}(X_{ij})^T\vbeta_{0j}\bigr) -
\rho_{\alpha}\bigl(Y_i^*\bigr) \bigr]
\\[-2pt]
&&\hspace*{-4pt}\qquad \defined I_1+I_2.
\end{eqnarray*}
Then $J_1\leq P(|I_1|\geq b_4n^{-2\tau}/4)+P(|I_2|\geq b_4n^{-2\tau}/4)$.
Note that $|W_i|\leq \break C |\bolds{\pi}(X_{ij})^T\vbeta_{0j} |$, for some positive constant $C$.
By the argument of Lemma~\ref{lemma1}, $\sup_{t}|f_j(t)-\bolds{\pi
}(X_{ij})^T\vbeta_{0j}|\leq c_2N^{-d}$. Thus,
$|W_i|$ are uniformly bounded by a constant $M$. Applying Bernstein's
inequality, there exists a positive
constant $b_5$ such that for all $n$ sufficiently large,
\begin{eqnarray*}
P \bigl( |I_1 |\geq b_4n^{-2\tau}/4 \bigr) \leq2
\exp \biggl(-\frac{b_4^2n^{1-4\tau}/16} {
2M^2+Mb_4n^{-2\tau}/3} \biggr) \leq2\exp \bigl(-b_5n^{1-4\tau}
\bigr).
\end{eqnarray*}
Furthermore, applying Lemma~\ref{lemma8},
\begin{eqnarray*}
I_2&=&n^{-2}\sum_{i=1}^n
\sum_{j=1}^n\delta_i\bigl[G
\bigl(Y_i^*\bigr)\bigr]^{-2}\xi \bigl(Y_j^*,
\delta_j,Y_i^*\bigr) \bigl[\rho_{\alpha}
\bigl(Y_i^*-\bolds{\pi}(X_{ij})^T
\vbeta_{0j}\bigr) -\rho_{\alpha}\bigl(Y_i^*\bigr)
\bigr]
\\
&&{}+n^{-1}\sum_{i=1}^n
\delta_iR_n\bigl(Y_i^*\bigr) \bigl[
\rho_{\alpha}\bigl(Y_i^*-\bolds{\pi}(X_{ij})^T
\vbeta_{0j}\bigr) -\rho_{\alpha}\bigl(Y_i^*\bigr)
\bigr]\defined I_{21}+I_{22},
\end{eqnarray*}
where $\xi$ and $R_n$ are defined in Lemma~\ref{lemma8}.
By Lemma~\ref{lemma8}, $I_{22}=\break O(n^{-3/4}(\log n)^{3/4})$ almost surely.
By assumptions $n^{-3/4}(\log n)^{3/4}=o(n^{-2\tau})$,
and noting that $\delta_iG^{-2}(Y_i^*)\xi(Y_j^*,\delta_j,Y_i^*)
[\rho_{\alpha}(Y_i^*-\bolds{\pi}(X_{ij})^T\vbeta_{0j})
-\rho_{\alpha}(Y_i^*) ]$ are independent bounded random variables,
we have for all $n$ sufficiently large,
\begin{eqnarray*}
&&\hspace*{-4pt}P \bigl( |I_2 |\geq b_4n^{-2\tau}/4 \bigr)
\\
&&\hspace*{-6pt}\qquad\leq P \Biggl( \frac{2}{n(n-1)}\sum_{i=1}^n
\sum_{j=1,j\neq i}^n\delta_iG^{-2}
\bigl(Y_i^*\bigr)\xi\bigl(Y_j^*,\delta_j,Y_i^*
\bigr) \\
&&\hspace*{-3pt}\hspace*{131pt}{}\times\bigl[\rho_{\alpha}\bigl(Y_i^*-\bolds{
\pi}(X_{ij})^T\vbeta_{0j}\bigr)
-\rho_{\alpha}\bigl(Y_i^*\bigr) \bigr]>b_4n^{-2\tau}/8
\Biggr)
\\
&&\hspace*{-6pt}\qquad \leq 2\exp \bigl(-b_6n^{1-4\tau} \bigr),
\end{eqnarray*}
where $b_6$ is a positive constant,
by Lemma~\ref{lemma9}.
Therefore, $J_1\leq4\exp (-c_{7}n^{1-4\tau} )$ where
$c_7=\min(b_5,b_6)$.

Next, we evaluate $J_2$.
%$D_n(\vbeta)=B_n(\vbeta)-B_n(\vbeta_{0j})-B(\vbeta)+B(\vbeta_{0j})$.
Let $V_i(\vbeta)=\rho_{\alpha}(Y_i^*-\bolds{\pi}(X_{ij})^T\vbeta
)-\rho_{\alpha}(Y_i^*-\bolds{\pi}(X_{ij})^T\vbeta_{0j})$
and let $Z_i=\delta_i[G(Y_i^*)]^{-1}V_i(\vbeta)$.
We have
\begin{eqnarray*}
J_2&\leq& P \Biggl(\sup_{\Vert \vbeta-\vbeta_{0j}\Vert \leq{CN^{1/2}n^{-\tau}}} \Biggl|n^{-1}\sum
_{i=1}^n \bigl[Z_i-E (Z_i
) \bigr] \Biggr| \\
&&\hspace*{14pt}{}+ \sup_{\Vert \vbeta-\vbeta_{0j}\Vert \leq{CN^{1/2}n^{-\tau}}} \Biggl|n^{-1}\sum_{i=1}^n
\delta_i\bigl[\bigl(\widehat {G}\bigl(Y_i^*\bigr)
\bigr)^{-1}-\bigl(G\bigl(Y_i^*\bigr)\bigr)^{-1}
\bigr]V_i(\vbeta)\Biggr |\\
&&\hspace*{243pt}{}\geq b_4n^{-2\tau}/2 \Biggr).
\end{eqnarray*}
Applying Knight's identity [(\citeyear{Knight}), page 758], we have
\begin{eqnarray*}
V_i(\vbeta)
&=& \bolds{\pi}(X_{ij})^T(\vbeta-\vbeta_{0j})
\bigl[I\bigl(Y_i^*-\bolds{\pi}(X_{ij})^T
\vbeta_{0j}\leq0\bigr)-\tau \bigr]
\\
&&{} + \int_0^{\bolds{\pi}(X_{ij})^T(\vbeta-\vbeta_{0j})} \bigl[I\bigl(Y_i^*-
\bolds{\pi}(X_{ij})^T\vbeta_{0j}\leq s\bigr) \\
&&\hspace*{86pt}{}-I
\bigl(Y_i^*-\bolds{\pi}(X_{ij})^T
\vbeta_{0j}\leq0\bigr) \bigr]\,ds.
\end{eqnarray*}
Thus,
\begin{eqnarray}
\label{wolf2}\qquad \sup_{\Vert \vbeta-\vbeta_{0j}\Vert \leq{CN^{1/2}n^{-\tau}}}\bigl|V_i(\vbeta )\bigr|&\leq& 2
\sup_{\Vert \vbeta-\vbeta_{0j}\Vert \leq{CN^{1/2}n^{-\tau}}} \bigl|\bolds{\pi}(X_{ij})^T(\vbeta-
\vbeta_{0j})\bigr |
\nonumber
\\[-8pt]
\\[-8pt]
\nonumber
&\leq& cNn^{-\tau}
\end{eqnarray}
for some $c>0$ because $\Vert B_k(\cdot)\Vert_{\infty}\leq1$.
Combining (\ref{wolf2}) with Lemma~\ref{lemma8}(3), we have
\begin{eqnarray*}
&&\sup_{\Vert \vbeta-\vbeta_{0j}\Vert \leq{CN^{1/2}n^{-\tau}}} \Biggl|n^{-1}\sum_{i=1}^n
\delta_i\bigl[\bigl(\widehat {G}\bigl(Y_i^*\bigr)
\bigr)^{-1}-\bigl(G\bigl(Y_i^*\bigr)\bigr)^{-1}
\bigr]V_i(\vbeta) \Biggr|
\\[-2pt]
&&\qquad = O\bigl(Nn^{-\tau-1/2}(\log n)^{1/2}\bigr)
\end{eqnarray*}
almost surely. Assume $N^2n^{2\tau-1}\log n=o(1)$, then for all $n$
sufficiently large,
\begin{eqnarray*}
J_2&\leq& P \Biggl(\sup_{\Vert \vbeta-\vbeta_{0j}\Vert \leq{CN^{1/2}n^{-\tau}}}\Biggl |n^{-1}\sum
_{i=1}^n \bigl[Z_i-E
(Z_i ) \bigr] \Biggr| \geq b_4n^{-2\tau}/4 \Biggr).
\end{eqnarray*}
We use Lemma~\ref{lemma7} to evaluate the above inequality. First note
that, (\ref{wolf2})
implies that
$\sup_{\Vert \vbeta-\vbeta_{0j}\Vert \leq CN^{1/2}n^{-\tau}}|Z_i|\leq
c^*Nn^{-\tau}$ for some positive constant $c^*$.
Next, let $e_1,\ldots,e_n$ be a Rademacher sequence (i.e., i.i.d.
sequence taking values of
$\pm$1 with probability $1/2$) independent of $Z_1,\ldots,Z_n$. We have
\begin{eqnarray*}
&&\E \Biggl\{\sup_{\Vert \vbeta-\vbeta_{0j}\Vert \leq{CN^{1/2}n^{-\tau}}} n^{-1} \Biggl|\sum
_{i=1}^n\bigl(Z_i-\E(Z_i)
\bigr) \Biggr| \Biggr\}
\\[-2pt]
&&\qquad\leq 2\E \Biggl\{\sup_{\Vert \vbeta-\vbeta_{0j}\Vert \leq
{CN^{1/2}n^{-\tau}}} n^{-1} \Biggl|\sum
_{i=1}^ne_iZ_i \Biggr| \Biggr\}
\\[-2pt]
&&\qquad \leq C\E \Biggl\{\sup_{\Vert \vbeta-\vbeta_{0j}\Vert \leq
{CN^{1/2}n^{-\tau}}} n^{-1} \Biggl|\sum
_{i=1}^n e_i\pi(X_{ij})^T(
\vbeta-\vbeta_{0j}) \Biggr| \Biggr\}
\\[-2pt]
&&\qquad \leq CN^{1/2}n^{-\tau}\E\Biggl\Vert n^{-1}\sum
_{i=1}^n e_i\pi(X_{ij})\Biggr\Vert
\leq CN^{1/2}n^{-\tau} \Biggl[\E\Biggl\Vert n^{-1}\sum
_{i=1}^n e_i
\pi(X_{ij})\Biggr\Vert^2 \Biggr]^{1/2}
\\[-2pt]
&&\qquad = CN^{1/2}n^{-\tau} \Biggl[n^{-2} \E \Biggl(\sum
_{i=1}^n e_i^2
\pi(X_{ij})^T\pi(X_{ij}) \Biggr)
\Biggr]^{1/2} \leq CN^{1/2}n^{-\tau-1/2}
\end{eqnarray*}
for some generic constant $C$ which may vary from line to line. In the above,
the first inequality applies the symmetrization theorem [Lemma 2.3.1,
\citet{Vaart}],
the second inequality applies the contraction theorem [Ledoux and
Talagrad, (\citeyear{Ledoux})] using the Lipschitz
property of the quantile objective function, and the last inequality uses
(\ref{bbound4}).
Now, we apply Lemma~\ref{lemma7} to evaluate $J_{2}$.
Let $Z=\sup_{\Vert \vbeta-\vbeta_{0j}\Vert \leq\Delta N^{1/2}n^{-\tau}} n^{-1}
|\sum_{i=1}^n(Z_i-\E(Z_i)) |$. In Lemma~\ref{lemma7}, we take
$t=b_4n^{-2\tau}/2-CN^{1/2}n^{-\tau-1/2}$ and
$L^2=4c^2N^{2}n^{-2\tau}$, which gives
\begin{eqnarray*}
J_{2}&=& P \bigl(Z\geq\E Z+\bigl(b_4n^{-2\tau}/4-\E
Z\bigr) \bigr)
\\[-2pt]
&\leq& P \bigl(Z\geq\E Z+\bigl(b_4n^{-2\tau}/4-CN^{1/2}n^{-\tau-1/2}
\bigr) \bigr)
\\[-2pt]
&\leq& \exp \biggl(-\frac{n(b_4n^{-2\tau}/4-CN^{1/2}n^{-\tau-1/2})^2} {
8c^2N^{2}n^{-2\tau}} \biggr)\leq\exp \bigl(-c_8N^{-2}n^{1-2\tau}
\bigr)
\end{eqnarray*}
for some positive constant $c_8$ and all $n$ sufficiently large.
\end{pf}

%s8.3 #&#
\subsection{\texorpdfstring{Proof of the Theorem \protect\ref{main2}}{Proof of the Theorem 4.1}}\label{sec8.3}
In this subsection, we establish the exponential tail probability
bounds for
$P ( |S_{jk} |\geq Cn^{-\tau}/8 )$, $k=1,\ldots,4$,
which lead to the
result of Theorem~\ref{main2}.

{\textit{An exponential tail probability bound for
$S_{j1}$}.}
Recall that
\begin{eqnarray*}
S_{j1}&=&\bigl(\widehat{\vbeta}_j^c-
\vbeta_{0j}\bigr)^T\bigl(\mathbb{P}_n\bolds {
\pi}(X_j)\bolds{\pi}(X_j)^T\bigr) \bigl(
\widehat {\vbeta}_j^c-\vbeta_{0j}\bigr)
\\
&\leq& \lambdamax\bigl(\mathbb{P}_n\bolds{\pi}(X_j)
\bolds{\pi }(X_j)^T\bigr)\bigl\Vert \widehat{
\vbeta}_j^c-\vbeta_{0j}\bigr\Vert^2.
\end{eqnarray*}

It follows from Lemmas {\ref{lemma3}} and~\ref{bear2} that for
some $C^*>0$,
%
%e8.8 #&#
\begin{eqnarray}
\label{dog} &&P \bigl(S_{j1}\geq Cn^{-\tau}/8 \bigr)
\nonumber
\\
&&\qquad\leq P \bigl(\lambdamax\bigl(\mathbb{P}_n\bolds{
\pi}(X_j)\bolds{\pi }(X_j)^T \bigr)
\geq(b_2+1)N^{-1}\bigr)
\nonumber
\\
&&\qquad\quad{}+P \bigl(\bigl\Vert \widehat{\vbeta}_j^c-
\vbeta_{0j}\bigr\Vert^2\geq (b_2+1)^{-1}CNn^{-\tau}/8
\bigr)
\nonumber
\\[-8pt]
\\[-8pt]
\nonumber
&&\qquad\leq 2N^2\exp \bigl(-c_4nN^{-3} \bigr)+P
\bigl(\bigl\Vert \widehat{\vbeta }_j^c-\vbeta_{0j}
\bigr\Vert >C^*N^{1/2}n^{-\tau/2} \bigr)
\nonumber
\\
&&\qquad\leq 2N^2\exp \bigl(-c_4nN^{-3} \bigr) +P
\bigl(\bigl\Vert \widehat{\vbeta}_j^c-\vbeta_{0j}
\bigr\Vert >C^*N^{1/2}n^{-\tau
} \bigr)
\nonumber
\\
&&\qquad\leq 2N^2\exp \bigl(-c_4nN^{-3} \bigr)+ 4
\exp \bigl(-c_7n^{1-4\tau} \bigr) +\exp \bigl(-c_{8}N^{-2}n^{1-2\tau}
\bigr).\nonumber
\end{eqnarray}

{\textit{An exponential tail probability bound for
$S_{j2}$}.}
We first establish an upper bound for $\Vert \vbeta_{0j}\Vert $.
By result (B.3) in the supplemental material [He, Wang and Hong (\citeyear{HWH})], $\E
[f_j(X_j)-f_{nj}(X_j)]^2\leq c_3N^{-2d}$, $\forall j$, for some $c_3>0$.
It follows that
\begin{eqnarray*}
\E \bigl[f_{nj}(X_j)^2 \bigr]&\leq& 2 \E
\bigl[f_{j}(X_j)^2 \bigr]+ 2\E \bigl[
\bigl(f_j(X_j)-f_{nj}(X_j)
\bigr)^2 \bigr]\\
&\leq& c_{9}+2c_3N^{-2d},
\end{eqnarray*}
for some positive constant $c_{9}$. Also note that
\begin{eqnarray*}
\E \bigl[f_{nj}(X_j)^2 \bigr]\geq\lambdamin
\bigl(\E\bolds{\pi }(X_j)\bolds{\pi}(X_j)^T
\bigr)\Vert \vbeta_{0j}\Vert^2\geq b_1N^{-1}
\Vert \vbeta_{0j}\Vert^2.
\end{eqnarray*}
This implies that
$
\Vert \vbeta_{0j}\Vert \leq c_{10}\sqrt{N}
$
for some positive constant $c_{10}$.

Since $|S_{j2}|\leq2\Vert \widehat{\vbeta}_j^c-\vbeta_{0j}\Vert \lambdamax
(\mathbb{P}_n\bolds{\pi}(X_j)\bolds{\pi}(X_j)^T)\Vert \vbeta_{0j}\Vert $,
we have
\begin{eqnarray*}
&&P \bigl(|S_{j2}|\geq Cn^{-\tau}/8 \bigr)
\\
&&\qquad\leq P \bigl(\bigl\Vert \widehat{\vbeta}_j^c-
\vbeta_{0j}\bigr\Vert \lambdamax \bigl(\mathbb{P}_n\bolds{
\pi}(X_j)\bolds{\pi}(X_j)^T\bigr) \geq
CN^{-1/2}n^{-\tau}/(16c_{10}) \bigr)
\\
&&\qquad\leq P \bigl(\lambdamax\bigl(\mathbb{P}_n\bolds{
\pi}(X_j)\bolds{\pi}(X_j)^T\bigr)
>(b_2+1)N^{-1} \bigr)
\\
&&\qquad\quad{}+P \bigl(\bigl\Vert \widehat{\vbeta}_j^c-
\vbeta_{0j}\bigr\Vert \geq (b_2+1)^{-1}CN^{1/2}n^{-\tau}/(16c_{10})
\bigr)
\\
&&\qquad\leq 2N^2\exp \bigl(-c_4nN^{-3} \bigr)+ 4
\exp \bigl(-c_7n^{1-4\tau} \bigr) +\exp \bigl(-c_{8}N^{-2}n^{1-2\tau}
\bigr).
\end{eqnarray*}

{\textit{An exponential tail probability bound for
$S_{j3}$}.}
We have
\begin{eqnarray*}
S_{j3}&=&\widehat{\vbeta}_j^c \bigl(
\mathbb{P}_n\bolds{\pi }(X_j)\bolds{\pi}(X_j)^T-
\E\bolds{\pi}(X_j)\bolds{\pi }(X_j)^T \bigr)
\vbeta_{0j}
\\
&=& \bigl(\widehat{\vbeta}_j^c-\vbeta_{0j}
\bigr) \bigl(\mathbb{P}_n\bolds{\pi }(X_j)\bolds{
\pi}(X_j)^T-\E\bolds{\pi}(X_j)\bolds{\pi
}(X_j)^T \bigr)\vbeta_{0j}
\\
&&{} +\vbeta_{0j}^T \bigl(\mathbb{P}_n\bolds{
\pi}(X_j)\bolds{\pi }(X_j)^T-\E\bolds{
\pi}(X_j)\bolds{\pi}(X_j)^T \bigr)
\vbeta_{0j}
\\
&\defined& S_{j31}+S_{j32}.
\end{eqnarray*}

Therefore,
\begin{eqnarray*}
&&P \bigl(|S_{j3}|\geq Cn^{-\tau}/8 \bigr)
\\
&&\qquad\leq P \bigl(|S_{j31}|\geq Cn^{-\tau}/16 \bigr) + P
\bigl(|S_{j32}|\geq Cn^{-\tau}/16 \bigr)
\\
&&\qquad\leq P \bigl(\bigl\Vert \widehat{\vbeta}_j^c-
\vbeta_{0j}\bigr\Vert \max\bigl( \bigl|\lambdamax(\vD_j) \bigr|, \bigl|
\lambdamin(\vD_j) \bigr|\bigr) \geq CN^{-1/2}n^{-\tau}/(16c_{10})
\bigr)
\\
&&\qquad\quad{}+P \bigl(\Vert \vbeta_{0j}\Vert^2\max\bigl( \bigl|\lambdamax(
\vD_j) \bigr|, \bigl|\lambdamin(\vD_j)\bigr |\bigr) \geq
Cn^{-\tau}/16 \bigr)
\\
&&\qquad\leq P \bigl(\max\bigl( \bigl|\lambdamax(\vD_j) \bigr|, \bigl|\lambdamin (
\vD_j) \bigr|\bigr) \geq N^{-1}/(16c_{10}) \bigr)\\
& &\qquad\quad{} +P \bigl(\bigl\Vert \widehat{\vbeta}_j^c-
\vbeta_{0j}\bigr\Vert \geq CN^{1/2}n^{-\tau
} \bigr)
\\
&&\qquad\quad{} +P \bigl(\max\bigl( \bigl|\lambdamax(\vD_j) \bigr|, \bigl|\lambdamin(
\vD_j)\bigr |\bigr) \geq CN^{-1}n^{-\tau}/
\bigl(16c_{10}^2\bigr) \bigr)
\\
&&\qquad\leq 2P \bigl(\max\bigl( \bigl|\lambdamax(\vD_j) \bigl|, \bigl|\lambdamin (
\vD_j) \bigr|\bigr) \geq C^*N^{-1}n^{-\tau}\bigr)
\\
&&\qquad\quad{}+P \bigl(\bigl\Vert \widehat{\vbeta}_j^c-
\vbeta_{0j}\bigr\Vert \geq CN^{1/2}n^{-\tau
} \bigr)
\\
&&\qquad\leq 2N^2\exp \bigl(-c_6N^{-3}n^{1-2\tau}
\bigr)+ 4\exp \bigl(-c_7n^{1-4\tau} \bigr) +\exp
\bigl(-c_{8}N^{-2}n^{1-2\tau} \bigr)
\end{eqnarray*}
for all $n$ sufficiently large, where the last inequality uses Lemmas~\ref{lemma3} and~\ref{bear2}.

{\textit{An exponential tail probability bound for $S_{j4}$}.}
\begin{eqnarray*}
S_{j4}&=& -2F_{\mathrm{KM},n}^{-1}(\alpha)n^{-1}
\sum_{i=1}^n \bigl[\bolds{\pi
}(X_j)^T\vbeta_{0j}-E\bolds{
\pi}(X_j)^T\vbeta_{0j} \bigr]
\\
&&{}-2F_{\mathrm{KM},n}^{-1}(\alpha)n^{-1}\sum
_{i=1}^n\bolds{\pi }(X_j)^T
\bigl(\widehat{\vbeta}_j^{c}-\vbeta_{0j}\bigr)
\\
&\defined& S_{j41}+S_{j42}.
\end{eqnarray*}
Note that $F_{\mathrm{KM},n}^{-1}(\alpha)$ is uniformly bounded for
$\beta_1\leq\alpha\leq\beta_2$ almost surely.
From the argument of Lemma~\ref{lemma1}, $E(\bolds{\pi
}(X_j)^T\vbeta_{0j})$ is uniformly bounded in $X_j$.
Applying Bernstein's inequality to $S_{j41}$, there exists a positive
constant $c_9$ such that
$P (|S_{j41}|>Cn^{-\tau}/16 )\leq\exp(-c_9n^{1-2\tau})$ for
all $n$ sufficiently large.
On the other hand, by the Cauchy--Schwarz inequality, for all $n$
sufficiently large,
\begin{eqnarray*}
&&P \bigl(|S_{j42}|>Cn^{-\tau}/16 \bigr)
\\
&&\qquad\leq P \Biggl(\Biggl|n^{-1}\sqrt{n} \Biggl[\sum
_{i=1}^n\bigl[\bolds{\pi }(X_j)^T
\bigl(\widehat{\vbeta}_j^{c}-\vbeta_{0j}\bigr)
\bigr]^2 \Biggr]^{1/2}\Biggr|>C^*n^{-\tau} \Biggr)
\\
&&\qquad\leq P \bigl( \bigl[ \bigl(\widehat{\vbeta}_j^{c}-
\vbeta_{0j}\bigr)^T \bigl(\mathbb {P}_n\bolds{
\pi}(X_j)\bolds{\pi}(X_j)^T \bigr) \bigl(
\widehat {\vbeta}_j^{c}-\vbeta_{0j} \bigr)
\bigr]^{1/2}>C^*n^{-\tau} \bigr)
\\
&&\qquad\leq P \bigl(\bigl\Vert \widehat{\vbeta}_j^{c}-
\vbeta_{0j}\bigr\Vert \lambda_{\max
}^{1/2}\bigl(
\mathbb{P}_n\bolds{\pi}(X_j)\bolds{\pi
}(X_j)^T\bigr)>C^*n^{-\tau} \bigr)
\\
&&\qquad\leq P \bigl(\lambdamax\bigl(\mathbb{P}_n\bolds{
\pi}(X_j)\bolds{\pi }(X_j)^T \bigr)
>(b_2+1)N^{-1}\bigr)
\\
&&\qquad\quad{}+P \bigl(\bigl\Vert \widehat{\vbeta}_j^c-
\vbeta_{0j}\bigr\Vert \geq C^*N^{1/2}n^{-\tau} \bigr)
\\
&&\qquad\leq 2N^2\exp \bigl(-c_6N^{-3}n^{1-2\tau}
\bigr)+ 4\exp \bigl(-c_7n^{1-4\tau} \bigr) +\exp
\bigl(-c_{8}N^{-2}n^{1-2\tau} \bigr)
\end{eqnarray*}
for all $n$ sufficiently large, where the last inequality uses Lemmas~\ref{lemma3} and~\ref{bear2} and $C^*$
denotes a generic positive constant which may vary from line to line.
Therefore, for all $n$ sufficiently large,
\begin{eqnarray*}
&&P \bigl(|S_{j4}|>Cn^{-\tau}/8 \bigr)
\\
&&\qquad\leq 2N^2\exp \bigl(-c_6N^{-3}n^{1-2\tau}
\bigr)+ 4\exp \bigl(-c_7n^{1-4\tau} \bigr) +\exp
\bigl(-c_{8}N^{-2}n^{1-2\tau} \bigr)
\\
&&\qquad\quad{}+\exp\bigl(-c_9n^{1-2\tau}\bigr)
\\
&&\qquad{}\leq 2N^2\exp \bigl(-c_6N^{-3}n^{1-2\tau}
\bigr)+ 5\exp \bigl(-c_7n^{1-4\tau} \bigr) +\exp
\bigl(-c_{8}N^{-2}n^{1-2\tau} \bigr).
\end{eqnarray*}

\begin{pf*}{Proof of Theorem~\ref{main2}}
\begin{longlist}[(1)]
\item[(1)] We have
\begin{eqnarray*}
&&P \Bigl(\max_{1\leq j\leq p} \bigl| \bigl\Vert \widehat {f}_{nj}^{c}
\bigr\Vert_n^2-\Vert f_{nj}\Vert^2\bigr | \geq
Cn^{-\tau} \Bigr)
\\
&&\qquad\leq p \bigl(4N^2\exp \bigl(-c_4N^{-3}n
\bigr)+ 17\exp \bigl(-c_7n^{1-4\tau} \bigr)+4\exp
\bigl(-c_{8}N^{-2}n^{1-2\tau} \bigr)
\\
&&\hspace*{168pt}\qquad\qquad{}+4N^2\exp \bigl(-c_6N^{-3}n^{1-2\tau}
\bigr) \bigr)
\\
&&\qquad\leq p \bigl(17\exp \bigl(-\delta_3n^{1-4\tau} \bigr)
+12N^2\exp \bigl(-\delta_{4}N^{-3}n^{1-2\tau}
\bigr) \bigr)
\end{eqnarray*}
for all $n$ sufficiently large, for some positive constants $\delta_3$
and $\delta_4$.

\item[(2)] The result follows by making use of the bound in (1)
and observing that
\begin{eqnarray*}
&&P \bigl(M_{\alpha}\subset\widehat{M}_{\alpha}^{c}
\bigr) \geq P \Bigl(\min_{j\in M_{\alpha}}\bigl\Vert \widehat{f}_{nj}^{c}
\bigr\Vert_n^2\geq \nu_n \Bigr)
\\
&&\qquad\geq P \Bigl(\min_{j\in M_{\alpha}}\Vert f_{nj}\Vert^2-
\max_{j\in
M_{\alpha}}\bigl | \bigl\Vert \widehat{f}_{nj}^{c}
\bigr\Vert_n^2-\Vert f_{nj}\Vert^2 \bigr|\geq
\nu_n \Bigr)
\\
&&\qquad =1- P \Bigl(\max_{j\in M_{\alpha}}\bigl | \bigl\Vert \widehat {f}_{nj}^{c}
\bigr\Vert_n^2-\Vert f_{nj}\Vert^2 \bigr|\geq
\min_{j\in M_{\alpha
}}\Vert f_{nj}\Vert^2-\nu_n
\Bigr)
\\
&&\qquad \geq1-P \Bigl(\max_{j\in M_{\alpha}} \bigl| \bigl\Vert \widehat {f}_{nj}^{c}
\bigr\Vert_n^2-\Vert f_{nj}\Vert^2 \bigr|\geq
c_1n^{-\tau}/16 \Bigr).
\end{eqnarray*}
\end{longlist}
\upqed\end{pf*}

\section*{Acknowledgements}
We are grateful to Professor B\"uhlmann, the Associate Editor and three referees for their
encouragement, insightful and constructive comments.

\begin{supplement}%[id=suppA]
\sname{Supplement A}\label{suppA}
\stitle{``Quantile-adaptive model-free variable screening for
high-dimensional heterogeneous data''}
\slink[doi]{10.1214/13-AOS1087SUPP} %[doi,text={...}] - jei reikia suskaldyti doi
\sdatatype{.pdf}
\sfilename{aos1087\_supp.pdf}
\sdescription{We provide additional technical details and numerical
examples in the
supplemental material.}
\end{supplement}

% imsref loaded by akundreckaite, 2013-01-31 15:51:10
% imsref loaded by akundreckaite, 2013-01-31 15:56:00
%


\printaddresses


\begin{thebibliography}{34}
% BibTex style file: ims.bst, 2013-01-28
% Default style options (sort=0,type=number).
% Used options (sort=1,type=nameyear).

%b1 #&#
\bibitem[\protect\citeauthoryear{Bair and Tibshirani}{2004}]{Bair}
%
\begin{barticle}[auto:STB|2013/01/29|08:09:18]
\bauthor{\bsnm{Bair},~\bfnm{E.}\binits{E.}} \AND
\bauthor{\bsnm{Tibshirani},~\bfnm{R.}\binits{R.}}
(\byear{2004}).
\btitle{Semi-supervised methods to predict patient survival from gene
expression data}.
\bjournal{PLoS Biol.}
\bvolume{2}
\bpages{511--522}.
\bptok{imsref}%
\end{barticle}
%
\endbibitem

%b2 #&#
\bibitem[\protect\citeauthoryear{Beran}{1981}]{Beran-1981}
%
\begin{bmisc}[auto:STB|2013/01/29|08:09:18]
\bauthor{\bsnm{Beran},~\bfnm{R.}\binits{R.}}
(\byear{1981}).
\bhowpublished{Nonparametric regression with randomly censored
survival data,
Technical report. Univ. California, Berkeley.}
\bptok{imsref}%
\end{bmisc}
%
\endbibitem

%b3 #&#
\bibitem[\protect\citeauthoryear{B{\"u}hlmann, Kalisch and
Maathuis}{2010}]{Buh-2010}
%
\begin{barticle}[mr]
\bauthor{\bsnm{B{\"u}hlmann},~\bfnm{P.}\binits{P.}},
\bauthor{\bsnm{Kalisch},~\bfnm{M.}\binits{M.}} \AND
\bauthor{\bsnm{Maathuis},~\bfnm{M.~H.}\binits{M.~H.}}
(\byear{2010}).
\btitle{Variable selection in high-dimensional linear models: Partially
faithful distributions and the {PC}-simple algorithm}.
\bjournal{Biometrika}
\bvolume{97}
\bpages{261--278}.
\bid{doi={10.1093/biomet/asq008}, issn={0006-3444}, mr={2650737}}
\bptok{imsref}%
\end{barticle}
%
\endbibitem

%b4 #&#
\bibitem[\protect\citeauthoryear{Fan, Feng and Wu}{2010}]{F12}
%
\begin{barticle}[auto:STB|2013/01/29|08:09:18]
\bauthor{\bsnm{Fan},~\bfnm{J.}\binits{J.}},
\bauthor{\bsnm{Feng},~\bfnm{Y.}\binits{Y.}} \AND
\bauthor{\bsnm{Wu},~\bfnm{Y.}\binits{Y.}}
(\byear{2010}).
\btitle{Ultrahigh dimensional variable selection for Cox's
proportional hazards
model}.
\bjournal{IMS Collections}
\bvolume{6}
\bpages{70--86}.
\bptok{imsref}%
\end{barticle}
%
\endbibitem

%b5 #&#
\bibitem[\protect\citeauthoryear{Fan, Feng and Song}{2011}]{F1}
%
\begin{barticle}[mr]
\bauthor{\bsnm{Fan},~\bfnm{Jianqing}\binits{J.}},
\bauthor{\bsnm{Feng},~\bfnm{Yang}\binits{Y.}} \AND
\bauthor{\bsnm{Song},~\bfnm{Rui}\binits{R.}}
(\byear{2011}).
\btitle{Nonparametric independence screening in sparse ultra-high-dimensional
additive models}.
\bjournal{J. Amer. Statist. Assoc.}
\bvolume{106}
\bpages{544--557}.
\bid{doi={10.1198/jasa.2011.tm09779}, issn={0162-1459}, mr={2847969}}
\bptok{imsref}%
\end{barticle}
%
\endbibitem

%b6 #&#
\bibitem[\protect\citeauthoryear{Fan and Li}{2001}]{FL}
%
\begin{barticle}[mr]
\bauthor{\bsnm{Fan},~\bfnm{Jianqing}\binits{J.}} \AND
\bauthor{\bsnm{Li},~\bfnm{Runze}\binits{R.}}
(\byear{2001}).
\btitle{Variable selection via nonconcave penalized likelihood and its oracle
properties}.
\bjournal{J. Amer. Statist. Assoc.}
\bvolume{96}
\bpages{1348--1360}.
\bid{doi={10.1198/016214501753382273}, issn={0162-1459}, mr={1946581}}
\bptok{imsref}%
\end{barticle}
%
\endbibitem

%b7 #&#
\bibitem[\protect\citeauthoryear{Fan and Lv}{2008}]{F2}
%
\begin{barticle}[auto:STB|2013/01/29|08:09:18]
\bauthor{\bsnm{Fan},~\bfnm{J.}\binits{J.}} \AND
\bauthor{\bsnm{Lv},~\bfnm{J.}\binits{J.}}
(\byear{2008}).
\btitle{Sure independence screening for ultra-high dimensional feature space
(with discussion)}.
\bjournal{J. Roy. Statist. Soc. Ser. B}
\bvolume{70}
\bpages{849--911}.
\bptok{imsref}%
\end{barticle}
%
\endbibitem

%b8 #&#
\bibitem[\protect\citeauthoryear{Fan, Samworth and Wu}{2009}]{F3}
%
\begin{barticle}[auto:STB|2013/01/29|08:09:18]
\bauthor{\bsnm{Fan},~\bfnm{J.}\binits{J.}},
\bauthor{\bsnm{Samworth},~\bfnm{R.}\binits{R.}} \AND
\bauthor{\bsnm{Wu},~\bfnm{Y.}\binits{Y.}}
(\byear{2009}).
\btitle{Ultrahigh dimensional variable selection: Beyond the linear model}.
\bjournal{J. Mach. Learn. Res.}
\bvolume{10}
\bpages{1829--1853}.
\bptok{imsref}%
\end{barticle}
%
\endbibitem

%b9 #&#
\bibitem[\protect\citeauthoryear{Fan and Song}{2010}]{F4}
%
\begin{barticle}[mr]
\bauthor{\bsnm{Fan},~\bfnm{Jianqing}\binits{J.}} \AND
\bauthor{\bsnm{Song},~\bfnm{Rui}\binits{R.}}
(\byear{2010}).
\btitle{Sure independence screening in generalized linear models with
{NP}-dimensionality}.
\bjournal{Ann. Statist.}
\bvolume{38}
\bpages{3567--3604}.
\bid{doi={10.1214/10-AOS798}, issn={0090-5364}, mr={2766861}}
\bptok{imsref}%
\end{barticle}
%
\endbibitem

%b10 #&#
\bibitem[\protect\citeauthoryear{Gonzalez-Manteiga and
Cadarso-Suarez}{1994}]{GMCS-1994}
%
\begin{barticle}[mr]
\bauthor{\bsnm{Gonzalez-Manteiga},~\bfnm{W.}\binits{W.}} \AND
\bauthor{\bsnm{Cadarso-Suarez},~\bfnm{C.}\binits{C.}}
(\byear{1994}).
\btitle{Asymptotic properties of a generalized {K}aplan--{M}eier
estimator with
some applications}.
\bjournal{J. Nonparametr. Stat.}
\bvolume{4}
\bpages{65--78}.
\bid{doi={10.1080/10485259408832601}, issn={1048-5252}, mr={1366364}}
\bptok{imsref}%
\end{barticle}
%
\endbibitem

%b11 #&#
\bibitem[\protect\citeauthoryear{Hall and Miller}{2009}]{Hall}
%
\begin{barticle}[mr]
\bauthor{\bsnm{Hall},~\bfnm{Peter}\binits{P.}} \AND
\bauthor{\bsnm{Miller},~\bfnm{Hugh}\binits{H.}}
(\byear{2009}).
\btitle{Using generalized correlation to effect variable selection in
very high
dimensional problems}.
\bjournal{J. Comput. Graph. Statist.}
\bvolume{18}
\bpages{533--550}.
\bid{doi={10.1198/jcgs.2009.08041}, issn={1061-8600}, mr={2751640}}
\bptok{imsref}%
\end{barticle}
%
\endbibitem

%b12 #&#
\bibitem[\protect\citeauthoryear{He and Shi}{1996}]{He2}
%
\begin{barticle}[mr]
\bauthor{\bsnm{He},~\bfnm{Xuming}\binits{X.}} \AND
\bauthor{\bsnm{Shi},~\bfnm{Peide}\binits{P.}}
(\byear{1996}).
\btitle{Bivariate tensor-product {$B$}-splines in a partly linear model}.
\bjournal{J. Multivariate Anal.}
\bvolume{58}
\bpages{162--181}.
\bid{doi={10.1006/jmva.1996.0045}, issn={0047-259X}, mr={1405586}}
\bptok{imsref}%
\end{barticle}
%
\endbibitem

%b13 #&#
\bibitem[\protect\citeauthoryear{He, Wang and Hong}{2013}]{HWH}
%
\begin{bmisc}[auto:STB|2013/01/29|08:09:18]
\bauthor{\bsnm{He},~\bfnm{X.}\binits{X.}},
\bauthor{\bsnm{Wang},~\bfnm{L.}\binits{L.}} \AND
\bauthor{\bsnm{Hong},~\bfnm{H.~G.}\binits{H.~G.}}
(\byear{2013}).
\bhowpublished{Supplement to ``Quantile-adaptive model-free variable screening
for high-dimensional heterogeneous data.''
DOI:\doiurl{10.1214/13-AOS1087SUPP}.}
\bptok{imsref}%
\end{bmisc}
%
\endbibitem

%b14 #&#
\bibitem[\protect\citeauthoryear{Hjort and Pollard}{1993}]{Hjort}
%
\begin{bmisc}[auto:STB|2013/01/29|08:09:18]
\bauthor{\bsnm{Hjort},~\bfnm{N.~L.}\binits{N.~L.}} \AND
\bauthor{\bsnm{Pollard},~\bfnm{D.}\binits{D.}}
(\byear{1993}).
\bhowpublished{Asymptotics for minimisers of convex processes. Technical
report, Dept. Statistics, Yale Univ., New Haven, CT. Available at
\url{http://citeseer.ist.psu.edu/hjort93asymptotics.html}.}
\bptok{imsref}%
\end{bmisc}
%
\endbibitem

%b15 #&#
\bibitem[\protect\citeauthoryear{Hoeffding}{1963}]{Hoeffding}
%
\begin{barticle}[mr]
\bauthor{\bsnm{Hoeffding},~\bfnm{Wassily}\binits{W.}}
(\byear{1963}).
\btitle{Probability inequalities for sums of bounded random variables}.
\bjournal{J. Amer. Statist. Assoc.}
\bvolume{58}
\bpages{13--30}.
\bid{issn={0162-1459}, mr={0144363}}
\bptok{imsref}%
\end{barticle}
%
\endbibitem

%b16 #&#
\bibitem[\protect\citeauthoryear{Knight}{1998}]{Knight}
%
\begin{barticle}[mr]
\bauthor{\bsnm{Knight},~\bfnm{Keith}\binits{K.}}
(\byear{1998}).
\btitle{Limiting distributions for {$L_1$} regression estimators under
general conditions}.
\bjournal{Ann. Statist.}
\bvolume{26}
\bpages{755--770}.
\bid{doi={10.1214/aos/1028144858}, issn={0090-5364}, mr={1626024}}
\bptok{imsref}%
\end{barticle}
%
\endbibitem

%b17 #&#
\bibitem[\protect\citeauthoryear{Koenker}{2005}]{K2}
%
\begin{bbook}[mr]
\bauthor{\bsnm{Koenker},~\bfnm{Roger}\binits{R.}}
(\byear{2005}).
\btitle{Quantile Regression}.
\bseries{Econometric Society Monographs}
\bvolume{38}.
\bpublisher{Cambridge Univ. Press}, \blocation{Cambridge}.
\bid{doi={10.1017/CBO9780511754098}, mr={2268657}}
\bptok{imsref}%
\end{bbook}
%
\endbibitem

%b18 #&#
\bibitem[\protect\citeauthoryear{Ledoux and Talagrand}{1991}]{Ledoux}
%
\begin{bbook}[mr]
\bauthor{\bsnm{Ledoux},~\bfnm{Michel}\binits{M.}} \AND
\bauthor{\bsnm{Talagrand},~\bfnm{Michel}\binits{M.}}
(\byear{1991}).
\btitle{Probability in {B}anach Spaces: Isoperimetry and Processes}.
\bseries{Ergebnisse der Mathematik und Ihrer Grenzgebiete (3) [Results in
Mathematics and Related Areas (3)]}
\bvolume{23}.
\bpublisher{Springer}, \blocation{Berlin}.
\bid{mr={1102015}}
\bptok{imsref}%
\end{bbook}
%
\endbibitem

%b19 #&#
\bibitem[\protect\citeauthoryear{Li and Luan}{2005}]{Li}
%
\begin{barticle}[pbm]
\bauthor{\bsnm{Li},~\bfnm{Hongzhe}\binits{H.}} \AND
\bauthor{\bsnm{Luan},~\bfnm{Yihui}\binits{Y.}}
(\byear{2005}).
\btitle{Boosting proportional hazards models using smoothing splines, with
applications to high-dimensional microarray data}.
\bjournal{Bioinformatics}
\bvolume{21}
\bpages{2403--2409}.
\bid{doi={10.1093/bioinformatics/bti324}, issn={1367-4803}, pii={bti324},
pmid={15713732}}
\bptok{imsref}%
\end{barticle}
%
\endbibitem

%b20 #&#
\bibitem[\protect\citeauthoryear{Li, Zhong and Zhu}{2012}]{LZZ}
%
\begin{barticle}[auto:STB|2013/01/29|08:09:18]
\bauthor{\bsnm{Li},~\bfnm{R.}\binits{R.}},
\bauthor{\bsnm{Zhong},~\bfnm{W.}\binits{W.}} \AND
\bauthor{\bsnm{Zhu},~\bfnm{L.}\binits{L.}}
(\byear{2012}).
\btitle{Feature screening via distance correlation learning}.
\bjournal{J. Amer. Statist. Assoc.}
\bvolume{107}
\bpages{1129--1139}.
\bptok{imsref}%
\end{barticle}
%
\endbibitem

%b21 #&#
\bibitem[\protect\citeauthoryear{Lo and Singh}{1986}]{Lo}
%
\begin{barticle}[mr]
\bauthor{\bsnm{Lo},~\bfnm{Shaw-Hwa}\binits{S.-H.}} \AND
\bauthor{\bsnm{Singh},~\bfnm{Kesar}\binits{K.}}
(\byear{1986}).
\btitle{The product-limit estimator and the bootstrap: Some asymptotic
representations}.
\bjournal{Probab. Theory Related Fields}
\bvolume{71}
\bpages{455--465}.
\bid{doi={10.1007/BF01000216}, issn={0178-8051}, mr={0824714}}
\bptok{imsref}%
\end{barticle}
%
\endbibitem

%b22 #&#
\bibitem[\protect\citeauthoryear{Massart}{2000}]{Massart}
%
\begin{barticle}[mr]
\bauthor{\bsnm{Massart},~\bfnm{Pascal}\binits{P.}}
(\byear{2000}).
\btitle{Some applications of concentration inequalities to statistics}.
\bjournal{Ann. Fac. Sci. Toulouse Math. (6)}
\bvolume{9}
\bpages{245--303}.
\bid{issn={0240-2963}, mr={1813803}}
\bptok{imsref}%
\end{barticle}
%
\endbibitem

%b23 #&#
\bibitem[\protect\citeauthoryear{McKeague, Subramanian and
Sun}{2001}]{McKeague-2001}
%
\begin{barticle}[mr]
\bauthor{\bsnm{McKeague},~\bfnm{Ian~W.}\binits{I.~W.}},
\bauthor{\bsnm{Subramanian},~\bfnm{Sundarraman}\binits{S.}} \AND
\bauthor{\bsnm{Sun},~\bfnm{Yanqing}\binits{Y.}}
(\byear{2001}).
\btitle{Median regression and the missing information principle}.
\bjournal{J. Nonparametr. Stat.}
\bvolume{13}
\bpages{709--727}.
\bid{doi={10.1080/10485250108832873}, issn={1048-5252}, mr={1931164}}
\bptok{imsref}%
\end{barticle}
%
\endbibitem

%b24 #&#
\bibitem[\protect\citeauthoryear{Peng and Huang}{2008}]{peng}
%
\begin{barticle}[mr]
\bauthor{\bsnm{Peng},~\bfnm{Limin}\binits{L.}} \AND
\bauthor{\bsnm{Huang},~\bfnm{Yijian}\binits{Y.}}
(\byear{2008}).
\btitle{Survival analysis with quantile regression models}.
\bjournal{J.~Amer. Statist. Assoc.}
\bvolume{103}
\bpages{637--649}.
\bid{doi={10.1198/016214508000000355}, issn={0162-1459}, mr={2435468}}
\bptok{imsref}%
\end{barticle}
%
\endbibitem

%b25 #&#
\bibitem[\protect\citeauthoryear{Portnoy}{2003}]{Portnoy-2003}
%
\begin{barticle}[mr]
\bauthor{\bsnm{Portnoy},~\bfnm{Stephen}\binits{S.}}
(\byear{2003}).
\btitle{Censored regression quantiles}.
\bjournal{J. Amer. Statist. Assoc.}
\bvolume{98}
\bpages{1001--1012}.
\bid{doi={10.1198/016214503000000954}, issn={0162-1459}, mr={2041488}}
\bptok{imsref}%
\end{barticle}
%
\endbibitem

%b26 #&#
\bibitem[\protect\citeauthoryear{Rosenwald et~al.}{2002}]{Ro}
%
\begin{barticle}[auto:STB|2013/01/29|08:09:18]
\bauthor{\bsnm{Rosenwald},~\bfnm{A.}\binits{A.}},
\bauthor{\bsnm{Wright},~\bfnm{G.}\binits{G.}},
\bauthor{\bsnm{Chan},~\bfnm{W.~C.}\binits{W.~C.}},
\bauthor{\bsnm{Connors},~\bfnm{J.~M.}\binits{J.~M.}},
\bauthor{\bsnm{Hermelink},~\bfnm{H.~K.}\binits{H.~K.}},
\bauthor{\bsnm{Smeland},~\bfnm{E.~B.}\binits{E.~B.}} \AND
\bauthor{\bsnm{Staudt},~\bfnm{L.~M.}\binits{L.~M.}}
(\byear{2002}).
\btitle{The use of molecular profiling to predict survival after chemotherapy
for diffuse large-B-cell lymphoma}.
\bjournal{The New England Journal of Medicine}
\bvolume{346}
\bpages{1937--1947}.
\bptok{imsref}%
\end{barticle}
%
\endbibitem

%b27 #&#
\bibitem[\protect\citeauthoryear{Stone}{1985}]{Stone}
%
\begin{barticle}[mr]
\bauthor{\bsnm{Stone},~\bfnm{Charles~J.}\binits{C.~J.}}
(\byear{1985}).
\btitle{Additive regression and other nonparametric models}.
\bjournal{Ann. Statist.}
\bvolume{13}
\bpages{689--705}.
\bid{doi={10.1214/aos/1176349548}, issn={0090-5364}, mr={0790566}}
\bptok{imsref}%
\end{barticle}
%
\endbibitem

%b28 #&#
\bibitem[\protect\citeauthoryear{van~der Vaart and Wellner}{1996}]{Vaart}
%
\begin{bbook}[mr]
\bauthor{\bparticle{van~der} \bsnm{Vaart},~\bfnm{Aad~W.}\binits
{A.~W.}} \AND
\bauthor{\bsnm{Wellner},~\bfnm{Jon~A.}\binits{J.~A.}}
(\byear{1996}).
\btitle{Weak Convergence and Empirical Processes}.
\bpublisher{Springer}, \blocation{New York}.
\bid{mr={1385671}}
\bptok{imsref}%
\end{bbook}
%
\endbibitem

%b29 #&#
\bibitem[\protect\citeauthoryear{Wang and Wang}{2009}]{wang}
%
\begin{barticle}[mr]
\bauthor{\bsnm{Wang},~\bfnm{Huixia~Judy}\binits{H.~J.}} \AND
\bauthor{\bsnm{Wang},~\bfnm{Lan}\binits{L.}}
(\byear{2009}).
\btitle{Locally weighted censored quantile regression}.
\bjournal{J.~Amer. Statist. Assoc.}
\bvolume{104}
\bpages{1117--1128}.
\bid{doi={10.1198/jasa.2009.tm08230}, issn={0162-1459}, mr={2562007}}
\bptok{imsref}%
\end{barticle}
%
\endbibitem

%b30 #&#
\bibitem[\protect\citeauthoryear{Ying, Jung and Wei}{1995}]{Ying-1995}
%
\begin{barticle}[mr]
\bauthor{\bsnm{Ying},~\bfnm{Z.}\binits{Z.}},
\bauthor{\bsnm{Jung},~\bfnm{S.~H.}\binits{S.~H.}} \AND
\bauthor{\bsnm{Wei},~\bfnm{L.~J.}\binits{L.~J.}}
(\byear{1995}).
\btitle{Survival analysis with median regression models}.
\bjournal{J.~Amer. Statist. Assoc.}
\bvolume{90}
\bpages{178--184}.
\bid{issn={0162-1459}, mr={1325125}}
\bptok{imsref}%
\end{barticle}
%
\endbibitem

%b31 #&#
\bibitem[\protect\citeauthoryear{Zhao and Li}{2012}]{Zhao}
%
\begin{barticle}[mr]
\bauthor{\bsnm{Zhao},~\bfnm{Sihai~Dave}\binits{S.~D.}} \AND
\bauthor{\bsnm{Li},~\bfnm{Yi}\binits{Y.}}
(\byear{2012}).
\btitle{Principled sure independence screening for {C}ox models with
ultra-high-dimensional covariates}.
\bjournal{J. Multivariate Anal.}
\bvolume{105}
\bpages{397--411}.
\bid{doi={10.1016/j.jmva.2011.08.002}, issn={0047-259X}, mr={2877525}}
\bptok{imsref}%
\end{barticle}
%
\endbibitem

%b32 #&#
\bibitem[\protect\citeauthoryear{Zhou, Shen and Wolfe}{1998}]{Zhou}
%
\begin{barticle}[mr]
\bauthor{\bsnm{Zhou},~\bfnm{S.}\binits{S.}},
\bauthor{\bsnm{Shen},~\bfnm{X.}\binits{X.}} \AND
\bauthor{\bsnm{Wolfe},~\bfnm{D.~A.}\binits{D.~A.}}
(\byear{1998}).
\btitle{Local asymptotics for regression splines and confidence regions}.
\bjournal{Ann. Statist.}
\bvolume{26}
\bpages{1760--1782}.
\bid{doi={10.1214/aos/1024691356}, issn={0090-5364}, mr={1673277}}
\bptok{imsref}%
\end{barticle}
%
\endbibitem

%b33 #&#
\bibitem[\protect\citeauthoryear{Zhu et~al.}{2011}]{Zhu}
%
\begin{barticle}[mr]
\bauthor{\bsnm{Zhu},~\bfnm{Li-Ping}\binits{L.-P.}},
\bauthor{\bsnm{Li},~\bfnm{Lexin}\binits{L.}},
\bauthor{\bsnm{Li},~\bfnm{Runze}\binits{R.}} \AND
\bauthor{\bsnm{Zhu},~\bfnm{Li-Xing}\binits{L.-X.}}
(\byear{2011}).
\btitle{Model-free feature screening for ultrahigh-dimensional data}.
\bjournal{J. Amer. Statist. Assoc.}
\bvolume{106}
\bpages{1464--1475}.
\bid{doi={10.1198/jasa.2011.tm10563}, issn={0162-1459}, mr={2896849}}
\bptok{imsref}%
\end{barticle}
%
\endbibitem

\end{thebibliography}
\end{document}